\documentclass[11pt]{amsart}

\usepackage{latexsym,amssymb,amsmath,amsopn,graphics,xy,epsfig,picture,epic}

\def\aa{{\mathcal A}}

\def\cc{{\mathcal C}}
\def\dd{{\mathcal D}}

\def\ff{{\mathcal F}}

\def\ss{{\mathcal S}}

\def\zzz{{\mathcal Z}}

\def\ffi{\varphi}
\def\eps{\varepsilon}
\def\dst{\displaystyle}

\DeclareMathOperator{\sinc}{sinc}

\renewcommand{\Im}{\mathrm{Im}\,}

\def\C{{\mathbb{C}}}

\def\N{{\mathbb{N}}}

\def\Q{{\mathbb{Q}}}
\def\R{{\mathbb{R}}}

\def\Z{{\mathbb{Z}}}

\newcommand{\norm}[1]{{\left\|{#1}\right\|}}
\newcommand{\ent}[1]{{\left[{#1}\right]}}
\newcommand{\abs}[1]{{\left|{#1}\right|}}
\newcommand{\scal}[1]{{\left\langle{#1}\right\rangle}}

\newtheorem{problem}{Problem}

\newtheorem{lemma}{Lemma}[section]
\newtheorem{proposition}[lemma]{Proposition}
\newtheorem{theorem}[lemma]{Theorem}

\theoremstyle{definition}
\newtheorem{definition}[lemma]{Definition}

\theoremstyle{remark}
\newtheorem{remark}[lemma]{Remark}

\begin{document}

\title[Uniqueness in phase retrival of FrFT's]{Uniqueness results for the phase retrieval problem
of fractional Fourier transforms of variable order}
\author{Philippe Jaming}
\address{Universit\'e d'Orl\'eans\\
Facult\'e des Sciences\\ 
MAPMO - F\'ed\'eration Denis Poisson\\ BP 6759\\ F 45067 Orl\'eans Cedex 2\\
France}
\curraddr{Institut de Math\'ematiques de Bordeaux UMR 5251,
Universit\'e Bordeaux 1, cours de la Lib\'eration, F 33405 Talence cedex, France}

\email{Philippe.Jaming@univ-orleans.fr}

\begin{abstract}
In this paper, we investigate the uniqueness of the phase retrieval problem
for the fractional Fourier transform (FrFT) of variable order.
This problem occurs naturally in optics and quantum physics. More precisely, we show that if
$u$ and $v$ are such that fractional Fourier transforms of order $\alpha$ have same modulus
$|\ff_\alpha u|=|\ff_\alpha v|$ for some set $\tau$ of $\alpha$'s, then $v$ is equal to $u$ up to
a constant phase factor. The set $\tau$ depends on some extra assumptions either on $u$ or on both $u$ and $v$. Cases considered here are $u$, $v$ of compact support, pulse trains, Hermite functions or linear combinations of translates and dilates of Gaussians. In this last case, the set $\tau$ may even be reduced to a single point ({\it i.e.} one fractional Fourier transform may suffice for uniqueness in the problem).
\end{abstract}



\maketitle

\section{Introduction}
Usually, when one measures a quantity, due to the nature of measurement equipment, noise, transmission in messy media... the phase of the quantity one wishes to know is lost.
In mathematical terms, one wants to know a quantity $\ffi(t)$ knowing only
$|\ffi(t)|$ for all $t\in\R^d$. Stated as this, the problem has too many
solutions to be useful and one tries to incorporate {\it a priori} knowledge on $\ffi$
to decrease the under-determination. Problems of that kind are called \emph{Phase Retrieval Problems} and arise in such diverse fields as microscopy ({\it see e.g.} \cite{37,47,61,83,118,119}), holography
\cite{42,115}, crystallography \cite{Mi,Ro}, neutron radiography \cite{3},
optical coherence tomography \cite{chandra}, optical design \cite{39}, radar signal processing \cite{Ja},
quantum mechanics \cite{Co,CH,Ja,Is,LR,Ma} to name a few. We refer to the books \cite{Hu,St}, the review articles \cite{KST,Mi,Fi,LBL} for descriptions of various instances of this problems, some solutions to it (both theoretical and numerical) and for further references.

The particular instance of the problem we are concerned with here deals with the so-called
\emph{Fractional Fourier Transform} (FrFT). Let us sketch a definition of this transform that is sufficient
for the needs of the introduction (a precise definition follows in Section \ref{sec:deffrac}).
First, for $u\in L^1(\R^d)\cap L^2(\R^d)$ we define the Fourier transform as
$$
\ff u(\xi)=\widehat{u}(\xi)=\int_{\R^d}u(t)e^{-2i\pi \scal{t,\xi}}\,\mbox{d}t,\quad\xi\in\R
$$
and then extend it to $L^2(\R^d)$ in the usual way. Here and throughout the paper $|\cdot|$ and $\scal{\cdot,\cdot}$ are respectively the standard Euclidean norm on $\R^d$ and the corresponding scalar product. The inverse Fourier transform is denoted
by $\ff^{-1}$. For $\alpha\in \R\setminus\pi\Z$, we define the fractional Fourier transform of
order $\alpha$ via
$$
\ff_\alpha u(\xi)=c_\alpha e^{-i\pi|\xi|^2\cot\alpha}\ff[e^{-i\pi|\cdot|^2\cot\alpha}u](\xi/\sin\alpha).
$$
where $c_\alpha$ is a normalisation constant. We define $\ff_0 u(\xi)=u(\xi)$, $\ff_\pi u(\xi)=u(-\xi)$.
Also note that, $\ff_{\pi/2}=\ff$, $\ff_{-\pi/2}=\ff^{-1}$ and that $\ff_\alpha\ff_\beta=\ff_{\alpha+\beta}$.

\medskip

This transform appears naturally in many instances including optics \cite{OZK},
quantum mechanics \cite{Ma,Na}, signal processing \cite{Al,OZK}...  
We will detail below several instances where the fractional Fourier transform occurs
and where one is further lead to the question of recovery of a function $u\in L^2(\R)$
from the phase-less measurements of several fractional Fourier transforms
$\{|\ff_{\alpha}u|\}_{\alpha\in \tau}$. 
More precisely, we deal with the following question:

\begin{problem}[Phase Retrieval Problem for the fractional Fourier transform]\label{pb:frft}\ \\
Let $u,v\in L^2(\R^d)$ and let $\tau\subset[0,\pi)$ be a set of indices (finite or not). Assume that $|\ff_\alpha v|=|\ff_\alpha u|$ for every $\alpha\in\tau$.
\begin{enumerate}
\renewcommand{\theenumi}{\roman{enumi}}
\item Does this imply that $v=cu$ for some constant $c\in\C$, $|c|=1$?

\item If we restrict $u\in\dd$ for some set $\dd\subset L^2(\R^d)$
do we then have $v=cu$ for some constant $c\in\C$, $|c|=1$?

\item If we further restrict both $u,v\in\dd$ for some set $\dd\subset L^2(\R^d)$
do we then have $v=cu$ for some constant $c\in\C$, $|c|=1$?
\end{enumerate}
In the first two cases we say that $u$ is uniquely determined (up to constant multiples or up to a 
constant phase factor) from
$\{|\ff_\alpha u|,\alpha\in\tau\}$. In the last case we say that  $u$ is uniquely determined (up to a
constant phase factor) from $\{|\ff_\alpha u|,\alpha\in\tau\}$ within the class $\dd$.
\end{problem}

The usual phase retrieval problem is the case $\tau=\{\pi/2\}$ ({\it i.e.} when $\ff_\alpha=\ff$
is the usual Fourier transform) within the class $\dd$ of compactly supported functions or distributions.
A part from this, the most famous problem of this sort is due to Pauli who asked whether
$|u|$ and $|\ff u|$ uniquely determines $u$ up to constant phase factors
{\it i.e.} here $\tau=\{0,\pi/2\}$.
Several counter-examples to this question have been constructed ({\it see e.g} \cite{CH,Is,Ja} and \cite{Co} for the state of the art on the problem).
In order to construct those counter-examples, it was shown that
$u$ is not uniquely determined from $|u|$ and $|\ff u|$ within the class $\{\sum_{\mathrm{finite}}a_i\gamma(x-x_i)\}$ where $\gamma$ is the standard Gaussian 
nor within the class
$\{\sum_{\mathrm{finite}}a_i\chi_{[0,1/2]}(x-i)\}$. We will show that for $\alpha$ well chosen,
$|u|$ and $|\ff_\alpha u|$ uniquely determines $u$ up to constant multiples within these classes.

Note that Reichenbach \cite{Re} conjectured that there is a unitary operator $U$ on $L^2(\R^d)$ such that
$|f|$, $|\hat f|$ and $|Uf|$ uniquely determine $f$ up to a constant phase factor. Our results thus show that
the fractional Fourier transform is a good candidate.

Further, we show that if the set of indices $\tau=[0,\pi)$, then $|\ff_\alpha u|$,
$\alpha\in\tau$ uniquely determines $u$ up to constant multiples and that the set $\tau$
can be reduced to a discrete set when $u$ and $v$ are compactly supported. Moreover, we provide theoretical
reconstruction formulae in this case. The numerical aspects of those formulae as well as
further algorithms are postponed to forthcoming work.

\medskip

This paper is organised as follows. We start with a section in which we present various instances of the 
fractional Fourier transform in physics. The following section is devoted to preliminaries on 
Fourier analysis and complex variables. We also solve the phase retrieval problem for one FrFT
phase-less measurement there in the class of compactly supported functions. 
We then devote Section \ref{sec:window}
 to two particular cases of phase-less measurements of windowed transforms (Fourier and wavelets). Section \ref{sec:reconstruction}
is devoted to the solution of Problem \ref{pb:frft} and is divided in several subsections dealing
with the three aspects of the problem.

\section{Physical models behind the Phase Retrieval Problem for the fractional Fourier transform}
\label{sec:physic}

In this section, we will present three instances of the phase retrieval problem for the fractional Fourier transform. The first and main one stems from optic measurements in the Fresnel domain. We will
here sketch how the fractional Fourier transform arises there. We then quickly present two more instances of our problem stemming from quantum physics.

\subsection{Diffraction phenomena in the optical far field}

An excellent derivation of how the fractional Fourier transform (or more precisely the Fresnel transform) occurs in optics can be found 
in \cite[Section 3.1]{LBL}. More detailed accounts can be found in \cite{BW,Go,So}.
We will here sketch the main features of the way the FrFT appears in optics, leaving aside
full mathematical rigour for which we refer to the above mentioned texts.
Roughly speaking, the fractional Fourier transforms
are adapted to the mathematical expression of the Fresnel
diffraction, just as the standard Fourier transform is adapted
to  Fraunhofer diffraction. This connection seems to have been made for the first time by
Pallat-Finet \cite{PF}. Let us now switch to the physical presentation.

\medskip

Light is an 
electromagnetic wave with coupled electric and magnetic fields travelling through space. In an homogeneous 
isotropic medium, like free space or a lens with constant refractive index, the electric and magnetic 
field vectors form a right-handed orthogonal triad with the direction of propagation. Disregarding 
polarisation, the field can be described by a scalar function $ \mathcal {U}(\mathbf{x};t)$ representing 
either the electric or the magnetic field amplitude. As the light used generally exhibits a strong 
monochromaticity, the time dependence of the field is a harmonic one and can thus be explicitly written 
as 
$$
\mathcal {U}(\mathbf{x};t)=\Re[U(\mathbf{x})e^{-i\omega t}].
$$
Here, $\omega$ denotes the angular frequency of the light, and the complex-valued amplitude or ``phasor''  
$U(\mathbf{x})$ depends on the spatial coordinates $\mathbf{x}= (x,y,z)\in\R^3$ only. If both quantities 
$\mathcal {U}(\mathbf{x};t)$ and $U(\mathbf{x})$ represent an optical wave, they must satisfy the
Helmholtz equation
\begin{equation}
\label{helmoltz}
\left(\Delta - \frac{1}{c^2} \frac{\partial^2}{\partial t^2}\right)\mathcal {U}(\mathbf{x};t)=0
\qquad{\text{and}}\qquad
\left( \Delta + k^2 \right) U(\mathbf{x}) = 0,
\end{equation}
where $\Delta=\frac{\partial^2}{\partial x^2}+\frac{\partial^2}{\partial y^2}+\frac{\partial^2}{\partial z^2}$ is the usual Laplacian. 
The Helmholtz equation directly follows from the Maxwell equations under the condition of a homogeneous 
medium and in absence of sources. The quantity $k$ is termed the \emph{wavenumber} or propagation
constant of the medium and is related to the light velocity $c$, the angular frequency $\omega$, and the 
vacuum wavelength $\lambda$ by
$$
k=\frac{\omega}{c} = \frac{2\pi \nu}{\lambda} \qquad\text{and}\qquad c = \frac{c_0}{\nu},
$$
whereby $\nu$ is the refractive index of the medium, {\it e.g.}, in vacuum $\nu_0$ = 1, and  $c_0$ is the
light velocity in vacuum. 

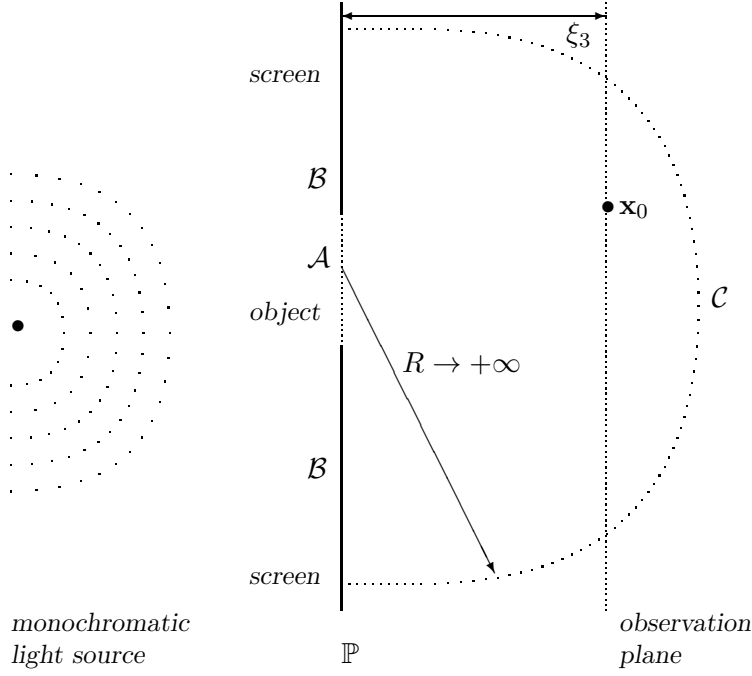
\begin{figure}
\begin{center}
\begin{picture}(350,250)
\put(0,125){$\bullet$}
\bezier{6}(0,105)(20,105)(20,125)
\bezier{6}(0,145)(20,145)(20,125)
\bezier{6}(0,95)(30,95)(30,125)
\bezier{6}(0,155)(30,155)(30,125)
\bezier{10}(0,85)(40,85)(40,125)
\bezier{10}(0,165)(40,165)(40,125)
\bezier{12}(0,75)(50,75)(50,125)
\bezier{12}(0,175)(50,175)(50,125)
\bezier{12}(0,65)(60,65)(60,125)
\bezier{12}(0,185)(60,185)(60,125)
\put(0,12){\sl\small monochromatic}
\put(0,0){\sl\small light source}
\put(112,150){$\mathcal{A}$}
\put(90,130){\sl\small object}
\put(112,70){$\mathcal{B}$}
\put(90,30){\sl\small screen}\put(90,220){\sl\small screen}
\put(112,180){$\mathcal{B}$}
\put(125,0){$\mathbb{P}$}
\drawline(125,20)(125,120)
\drawline(125,170)(125,250)
\dottedline{2}(125,120)(125,170)
\dottedline{2}(225,20)(225,250)
\put(230,12){\sl\small observation}
\put(230,0){\sl\small plane}
\put(175,245){\vector(-1,0){50}}
\put(175,245){\vector(1,0){50}}
\put(125,150){\vector(1,-2){58}}
\put(148,110){$R\to+\infty$}
\put(210,235){$\xi_3$}
\put(223,170){$\bullet$}
\put(230,170){$\mathbf{x}_0$}
\dottedline{3}(125,30)(155,30)
\dottedline{3}(125,240)(155,240)
\bezier{40}(155,30)(260,30)(260,135)
\bezier{40}(155,240)(260,240)(260,135)
\put(265,135){$\mathcal{C}$}
\end{picture}
\caption{Diffraction by a planar object}
\end{center}
\end{figure}

We will now study diffraction effects occurring from a
planar object placed in an opening of a plane screen $\mathbb{P}$.
as illustrated in the previous picture. In order to do so, we will appeal to Green's theorem: 
\begin{equation}
\label{eq:green}
-\int_{\mathcal{S}} \left(\tilde U(\mathbf{x})\frac{\partial U}{\partial\vec{n}}(\mathbf{x})
-U(\mathbf{x})\frac{\partial\tilde U}{\partial\vec{n}}(\mathbf{x})\right)\,\mbox{d}\mathbf{x}=
\iint_\Omega \tilde U(\mathbf{x})\Delta U(\mathbf{x})-U(\mathbf{x})\Delta\tilde U(\mathbf{x})\,\mbox{d}\mathbf{x}
\end{equation}
where $\mathcal{S}$ is the boundary of a domain $\Omega$ and is smooth closed and orientable and where
$\vec{n}$ is the \emph{inner} normal to that boundary.

The integration surface $\mathcal{S}$ is segmented into three disjoint parts, 
{\it i.e.}, $\mathcal{S}=\mathcal{A}\cup\mathcal {B}\cup\mathcal {C}$. The boundary $\mathcal{A}$ is 
chosen across the object location in the screen plane $\mathbb{P}$, $\mathcal {B}$ is the opaque part of the screen $\mathcal {B}=\mathbb{P}\setminus\mathcal{A}$, 
and $\mathcal{C}$ is the boundary of a half-sphere containing the observation point $\mathbf{x}_0$. 

Let us now  consider the unit-amplitude spherical wave expanding about the observation point $\xi$:
$$
G_0(\mathbf{x};\xi) = \frac{e^{jk_0\norm{\xi-\mathbf{x}}}}{\norm{\xi-\mathbf{x}}}
$$
and note that
$$
(\Delta +k_0^2)G_0(\mathbf{x};\xi)=4\pi\delta_{\xi-\mathbf{x}}.
$$
Let $\mathbf{x}_1$ be the symmetric of $\mathbf{x}_0$ with respect to $\mathbb{P}$
so that $\norm{\mathbf{x}_0-\mathbf{x}}=\norm{\mathbf{x}_1-\mathbf{x}}$ for all
$\mathbf{x}\in\mathbb{P}$ and define
$$
\tilde U(\mathbf{x})=G_0(\mathbf{x};\mathbf{x}_0)-G_0(\mathbf{x};\mathbf{x}_1).
$$
Then $\tilde U$ satisfies the following conditions
\begin{eqnarray}
(\Delta+k^2)\tilde U(\mathbf{x})=4\pi\bigl(\delta_{\mathbf{x}-\mathbf{x}_0}-\delta_{\mathbf{x}-\mathbf{x}_1}\bigr)&&\mbox{in }\Omega\\
\tilde U(\mathbf{x})=0&&\mbox{on }\mathbb{P}\\
\norm{\mathbf{x}-\mathbf{x}_0}\left(\frac{\partial\tilde U}{\partial\vec{n}}-ik\tilde U\right)\to 0&&\mbox{as }\norm{\mathbf{x}-\mathbf{x}_0}\to\infty.
\end{eqnarray}
From this, and the fact that $U$ also satisfies the Helmholtz equation, we derive from \eqref{eq:green}
the so-called integral theorem of Helmholtz and Kirchhoff ({\it see e.g.} \cite[p. 377]{BW}):
\begin{equation}
\label{eq:Helm-Kirch}
U(\mathbf{x}_0) = -\frac{1}{4\pi} \iint\limits_\mathcal{S}
\left( \frac{\partial U(\mathbf{x})}{\partial \vec{n}}\tilde U(\mathbf{x})
- U(\mathbf{x}) \frac{\partial \tilde U(\mathbf{x})}{\partial \vec{n}} \right)\,\mathrm{d}\mathbf{x}.
\end{equation}
Now, while $\tilde U$ is identically $0$ on the plane $\mathbb{P}$, its normal derivative is not.
The solution of the diffraction problem is found by specifying Kirchhoff boundary conditions on 
$\mathcal{A}$ and $\mathcal{B}$, and the Sommerfeld radiation condition on $\mathcal{C}$:
\begin{eqnarray*}
\mbox{ on }\mathcal{A}:&&U(\mathbf{x})=U_s(\mathbf{x})\quad\mbox{and}\quad\frac{\partial U}{\partial\vec{n}}(\mathbf{x})=\frac{\partial U_s}{\partial\vec{n}}(\mathbf{x})\\
\mbox{ on }\mathcal{B}:&&U(\mathbf{x})=0\quad\mbox{and}\quad\frac{\partial U}{\partial\vec{n}}(\mathbf{x})=0\\
\mbox{ on }\mathcal{C}:&&\lim_{\norm{x}\to+\infty}\norm{x}\left(\frac{\partial U}{\partial\vec{n}}(\mathbf{x})-ik U(\mathbf{x})\right)=0
\end{eqnarray*}
whereby $U_s(\mathbf{x})$ is the field incident on the screen. Let us comment on these assumptions.
The Kirchhoff boundary conditions are an idealization to the real field distribution on the screen by assuming that the field and its derivative across $\mathcal {A}$ is exactly the same as they would be in the absence of the screen. In the geometrical shadow across $\mathcal {B}$ the field is simply set to zero. Although these assumptions are reasonable, they fail in the immediate neighborhood of the rim of the opening.
The Sommerfeld radiation condition on $\mathcal {C}$ is satisfied, if the disturbance $U(\mathbf{x})$ vanishes at least as fast as a spherical wave. 

Integrating these boundary conditions into the integral \eqref{eq:Helm-Kirch} and letting the radius of the half-sphere $\mathcal{C}$ go to infinity, \eqref{eq:Helm-Kirch} yields
\begin{equation}
\label{eq:Helm-Kirch2}
U(\mathbf{x}_0) = -\frac{1}{4\pi} \iint\limits_\mathcal{A} U_s(\mathbf{x}) \frac{\partial G(\mathbf{x})}{\partial \vec{n}}\,\mathrm{d}\mathbf{x}.
\end{equation}

Now, a simple computation shows that
\begin{eqnarray*}
\frac{\partial\tilde U}{\partial\vec{n}}&=&2\frac{\exp(ik\norm{\mathbf{x}-\mathbf{x}_0})}{\norm{\mathbf{x}-\mathbf{x}_0}}
\left(ik-\frac{1}{\norm{\mathbf{x}-\mathbf{x}_0}}\right)
\frac{\scal{\vec{n},\mathbf{x}-\mathbf{x}_0}}{\norm{\mathbf{x}-\mathbf{x}_0}}\\
&\simeq&2ki\frac{\exp(ik\norm{\mathbf{x}-\mathbf{x}_0})}{\norm{\mathbf{x}-\mathbf{x}_0}}
\frac{\scal{\vec{n},\mathbf{x}-\mathbf{x}_0}}{\norm{\mathbf{x}-\mathbf{x}_0}}
\end{eqnarray*}
if $\norm{\mathbf{x}-\mathbf{x}_0}\gg\lambda$.

At this stage, it is useful to introduce the \emph{small angle} approximation wherein
the angle between $\vec{n}$ and $\mathbf{x}-\mathbf{x}_0$ is small, that is
$\dst \frac{\scal{\vec{n},\mathbf{x}-\mathbf{x}_0}}{\norm{\mathbf{x}-\mathbf{x}_0}}\simeq1$. For this, we establish reference coordinates $(x_1,x_2,x_3)$ relative to the plane $\mathbb{P}$, centered on the region $\aa$. We take the $x_3$ axis to be perpendicular to the plane $\mathbb{P}$. Further, we assume that the measurement plane in which $\mathbf{x}_0:=(\xi_1,\xi_2,\xi_3)$ lies is far from the screen plane $\mathbb{P}$. In other words, if $\mathbf{x}=(x_1,x_2,x_3)\in\mathbb{P}$, then $x_3=0$ and
$\norm{(x_1-\xi_1,x_2-\xi_2,0)}\ll\xi_3$. But then
\begin{eqnarray*}
\norm{\mathbf{x}-\mathbf{x}_0}&=&\xi_3\left(1+\frac{(\xi_1-x_1)^2}{\xi_3^2}+\frac{(\xi_2-x_2)^2}{\xi_3^2}\right)^{1/2}\\
&\simeq&\xi_3\left(1+\frac{(\xi_1-x_1)^2}{2\xi_3^2}+\frac{(\xi_2-x_2)^2}{2\xi_3^2}\right).
\end{eqnarray*}
Finally, under the various approximations made so far,
$$
\frac{\partial\tilde U}{\partial\vec{n}}\simeq\frac{\exp(ik\xi_3)}{i\lambda\xi_3}\exp\left(\frac{ik}{2\xi_3}\bigl((\xi_1-x_1)^2+(\xi_2-x_2)^2\bigr)\right).
$$

It follows that the field $U$ at an observation point $\mathbf{x}_0$ far enough from the screen is expressed in terms of the field $U_s$ on the screen as
$$
U(\mathbf{x}_0) = -\frac{\exp(ik\xi_3)}{4i\pi\lambda\xi_3} \iint\limits_\mathcal{A} U_s(\mathbf{x}) 
\exp\left(\frac{ik}{2\xi_3}\bigl((\xi_1-x_1)^2+(\xi_2-x_2)^2\bigr)\right)\,\mathrm{d}\mathbf{x}.
$$
A simple computation then shows that
$$
U(\mathbf{x}_0) =\frac{\exp(ik\xi_3)}{2i\lambda(2\pi\xi_3+ik)}e^{-i\pi(\xi_1^2+\xi_2^2)\sin^2\alpha}\ff_\alpha[U_s](\xi_1,\xi_2)
$$
with $\dst\cot\alpha=-\frac{k}{2\pi\xi_3}$.

Further occurrences of the fractional Fourier transform in optics may be found {\it e.g.} in \cite{Lo,MO1,MO2,OM,OZK,Z} and in the references therein.

\medskip

As optical measurement devices are not sensitive to phase, we are lead to 
the question of reconstruction of $U_s$ from the measured quantity
$|U(\mathbf{x}_0)|$; or equivalently
$$
\frac{1}{2\lambda\sqrt{4\pi^2\xi_3^2+k^2}}|\ff_{-\arg\cot k/(2\pi\xi_3)}[U_s](\xi_1,\xi_2)|.
$$
This leads us naturally to Problem \ref{pb:frft}, which we rephrase
as follows: {\it can one reconstruct a function $f$ from the phase-less measurements
$|\ff_{\alpha_1}f|$,... $|\ff_{\alpha_N}f|$ with $0<\alpha_1<\cdots<\alpha_N$.}

The problem of a single measurement seems not to have attracted any specific attention.
Indeed, the reader will easily check that this reduces to the case $\alpha_1=\pi/2$ {\it i.e.} the usual 
phase retrieval problem for the Fourier transform, thus a full description of the solutions can be given in terms of the so-called \emph{zero-flipping} phenomena. For sake of completeness, we have included a detailed proof in Section \ref{sec:phafr} below.

Generalisation to the case in which the intensity of an object and that of its
fractional Fourier transform, or the intensity of any
two of its fractional Fourier transforms, is known has
been addressed in \cite{15,16,17,18} to name a few.
However, the problem of retrieval from more than two Fresnel
transforms does not seem to have received anywhere
near the attention received by the two-intensity
problem and we are only aware of the papers \cite{48,49,50,51,52,EAOB}.
The problem
of recovery from multiple-intensity observations
related through general linear transformations is discussed
in \cite{ISV}.

However, from the mathematical point of view, those papers do not present a
rigorous proof that the object can be reconstructed from phase-less measurements,
but mainly show, using numerical exploration, that an adaptation of the Gerchberg-Saxton
algorithm will provide a reasonable solution. It is our aim here to provide a
complete mathematical proof of facts of that nature.

\subsection{Quantum Tomography}
Recall that the Wigner transform is defined for $u\in L^2(\R)$ by
$$
W(u)(\xi,x)=\int_{-\infty}^{+\infty}u(x+t/2)\overline{u(x-t/2)}e^{2i\pi\xi t}\,\mbox{d}t.
$$
An important property of the Wigner transform is that its marginals are given by:
$$
\int_{-\infty}^{+\infty}W(u)(\xi,x)\,\mbox{d}x=|\widehat{u}(\xi)|^2
\quad\mbox{and}
\int_{-\infty}^{+\infty}W(u)(\xi,x)\,\mbox{d}\xi=|u(x)|^2.
$$
Because of this property, this transform was proposed as a joint time-frequency probability distribution with marginals the distributions of position and velocity probability. Unfortunately $W(u)$ always takes negative values, unless $u$ is a Gaussian, so that this interpretation has some defect.

Using the link between the Wigner distribution and the ambiguity function and the fact that the ambiguity function of a fractional Fourier transform ({\it see} Section \ref{sec:amb} for details), one can check that marginals in other directions are given by the modulus of a fractional Fourier transform. More precisely:
$$
\int_{-\infty}^{+\infty}W(u)(t\cos\alpha-r\sin\alpha,t\sin\alpha+r\cos\alpha)\,\mbox{d}r
=|\ff_\alpha u(t)|^2.
$$
In particular, if we introduce the Radon transform of a function $F\in L^1(\R^2)\cap L^2(\R^2)$ on $\R^2$ as
$$
\mathcal{R}[F](\alpha,t)=\int_{x\,:\scal{x,(\cos\alpha,\sin\alpha)}=t}F(x)\,\mbox{d}x
=\int_\R F(t\cos\alpha-r\sin\alpha,t\sin\alpha+r\cos\alpha)\,\mbox{d}r
$$
(note that this function is only defined for almost all $t$). Then the previous relation reads
\begin{equation}
\label{eq:wignerradon}
\mathcal{R}[W(u)](\alpha,t)=|\ff_\alpha u(t)|^2.
\end{equation}
Because of this property, the problem of reconstructing $u$ from $|\ff_\alpha u(t)|^2$
for various $\alpha$'s is equivalent to the quantum state tomography problem
({\it see} \cite{LR} and references therein for more on this problem and \cite{Ma}
for the link with the fractional Fourier transform):

\begin{problem}[Quantum State Tomography]\ \\
Let $u,v\in L^2(\R)$ and let $\tau\subset[-\pi/2,\pi/2]$ be a set of angles.
Assume that
$$
\mathcal{R}[W(u)](\alpha,t)=\mathcal{R}[W(v)](\alpha,t)\quad\mathrm{for\ all}\ \alpha\in\tau\mathrm{\ and\ all\ }t\in\R.
$$
\begin{enumerate}
\renewcommand{\theenumi}{\roman{enumi}}
\item Does this imply that $v=cu$ for some constant $c\in\C$, $|c|=1$?

\item If we restrict $u\in\dd$ for some set $\dd\subset L^2(\R)$
do we then have $v=cu$ for some constant $c\in\C$, $|c|=1$?

\item If we further restrict both $u,v\in\dd$ for some set $\dd\subset L^2(\R^d)$
do we then have $v=cu$  for some constant $c\in\C$, $|c|=1$?
\end{enumerate}
\end{problem}

It follows that uniqueness and reconstruction results from partial Radon data may be transferred
into similar results for the Phase retrieval problem.
However, such techniques may not always be applicable. The main issue here is to have sufficient
``concentration'' of $W(u)$ near the origin so that it may be reconstructed via inversion of its Radon transform. Such concentration is unfortunately not always available, for instance, $W(u)$ is never compactly supported, unless $u=0$, so that most uniqueness results for the Radon transform
do not transfer directly into uniqueness results for the Phase Retrieval Problem
for fractional Fourier transforms.

A further issue is that of practical recovery by appealing to reconstruction algorithms
for the Radon transform such as filtered back projection ({\it see e.g.} \cite[Section V]{Nat}
or \cite[Section 3.6]{RK}). Recall that this algorithm allows to reconstruct a function
$W_{FBP}(u)$ from the data $\mathcal{R}[W(u)](\alpha_k,t_\ell)$, $\alpha_k=k\pi/N$, $t_\ell=\ell/N$
for some $N$ large enough. Note that, according to \eqref{eq:wignerradon}, this amounts
for reconstruct $W(u)$ from the data $|\ff_{\alpha_k} u(t_\ell)|$. This function $W_{FBP}(u)$ is then
an approximation of a smoothed version $\Phi*W(u)$ (where $*$ is the convolution over $\R^d$) of $W(u)$,
up to an error term of the form
$$
\lesssim \sup_{\theta\in[0,2\pi)}\int_{r\geq R}r|\ff[W(u)](r\cos\theta,r\sin\theta)|\mbox{d}r
$$
where $R$ depends on $N$, the ``sampling scale''. Of course, there is still the issue of reconstruction $u$ accurately from $\Phi*W(u)$. Note also that a large number of ``angular'' measurements $\alpha_k$
is needed here, an issue that may imply that such an algorithm may not be used in practice.
We postpone a detailed study of this algorithm to a forthcoming paper.

\subsection{Free Shr\"odinger Equation}
The formulation of the problem in this section corresponds to a question asked to the author by L. Vega \cite{Ve}. It can also be found rather implicitly in \cite{LR}.

Let us recall that the solution of the Free Shr\"odinger Equation
\begin{equation}
\label{eq:shr}
\begin{cases}\dst
i\partial_t u+\frac{1}{4\pi}\Delta_x^2 u=0\\
u(x,0)=u_0(x)\\
\end{cases}
\end{equation}
with initial data $u_0\in L^2(\R)$ has solution
$$
u(x,t)=\int_{\R}e^{-i\pi|\xi|^2t+2i\pi\scal{x,\xi}}\widehat{u_0}(\xi)\,\mbox{d}\xi
=\ff^{-1}\bigl[e^{-i\pi\xi^2t}\widehat{u_0}\bigr]\check\ (x).
$$
A straightforward computation then shows that, for $\alpha\in(-\pi/2,\pi/2)$,
$$
\ff_\alpha u_0(\xi)=
\left(\frac{ie^{i\alpha/2}}{\sqrt{|\cos\alpha|}}\right)^de^{-i\pi|\xi|^2\cot\alpha}u(\xi/\cos\alpha,\cot\alpha).
$$
We may thus rephrase the Phase Retrieval Problem as follows:

\begin{problem}[Phase Retrieval Problem for the Free Shr\"odinger Equation]\ \\
Let $u_0,v_0\in L^2(\R^d)$ and let $u$ and $v$ be the solutions of the
Free Shr\"odinger equation \eqref{eq:shr} with initial value $u_0$ and $v_0$.
Let $\tau=\{t_i\}_{i\in I}\subset[0,+\infty)$ be a set of times (finite or not)
at which one measures $|u|$ and $|v|$. Assume that
$$
|u(x,t_i)|=|v(x,t_i)|\quad\mathrm{for all}\ i\in I.
$$
\begin{enumerate}
\renewcommand{\theenumi}{\roman{enumi}}
\item Does this imply that $u=cv$ for some constant $c\in\C$, $|c|=1$?

\item If we restrict $u\in\dd$ for some set $\dd\subset L^2(\R^d)$
do we then have $u=cv$ for some constant $c\in\C$, $|c|=1$?

\item If we further restrict both $u,v\in\dd$ for some set $\dd\subset L^2(\R^d)$
do we then have $u=cv$ for some constant $c\in\C$, $|c|=1$?
\end{enumerate}
\end{problem}

Of course, a positive answer to any of these questions implies a positive answer to the subsequent ones.

Note that Pauli's problem whether an $u_0\in L^2(\R^d)$ is uniquely determined (up to constant multiples)
by $|u_0|$ and $|\hat u_0|$ may then be seen as a particular case of the above problem with $\tau=\{0,\infty\}$.

\section{Preliminaries}

\subsection{Reconstruction of an entire function of finite order from
its modulus on two lines}

In this section we gather some information about entire functions of finite order. We then prove
that an entire function of finite type is determined up to a constant phase factor by its modulus on two well chosen lines through the origin. Although this result is not of fundamental importance
in this paper and will only be applied in Section \ref{sec:window}, it is at the heart of the philosophy 
of the paper. Indeed, the fractional Fourier transform $\ff_\alpha u$ may be interpreted as the restriction of a function $u$ to the line $\ell_\alpha=\{(r\cos\alpha,r\sin\alpha)\in\R^2,\ r\in\R\}$ in the ``time-frequency plane'' ({\it see} Section \ref{sec:amb}). We thus want to reconstruct $u$ from its modulus on a few lines in the plane.
Of course, in general $\{\ff_\alpha u,\alpha\in\tau\}$ is not the restriction
of a single entire function $U$ to the lines $\{\ell_\alpha,\alpha\in\tau\}$ so that
the result presented here does not directly solve Problem \ref{pb:frft}.

\subsubsection{Preliminaries on entire functions of finite order}
Results in this section may be found in most books on one-variable complex analysis.

\begin{definition}
Let $f$ be an entire function ({\it i.e.} a function that is holomorphic over $\C$).
For $r>0$ define $\dst M_f(r)=\max_{|z|=r}|f(z)|$. The order $\rho$ of $f$ is defined by
$$
\rho=\limsup_{r\to+\infty}\frac{\log\log M_f(r)}{\log r}.
$$
\end{definition}
In other words, $f$ is of finite order $\rho$ if $\rho$ is the infimum of the non-negative numbers $\lambda$ such that
$M_f(r)\leq e^{r^{\lambda+\eps}}$ for any $\eps>0$ and $r$ sufficiently large.

Now, for $k\in\N$ and $\zeta\in\C$, let us denote by
$$
E_k(\zeta,z)=\begin{cases}\left(1-\frac{z}{\zeta}\right)&\mbox{if }k=0\\
\left(1-\frac{z}{\zeta}\right)\exp\left(\frac{z}{\zeta}+\frac{1}{2}\frac{z^2}{\zeta^2}\cdots+
\frac{1}{k}\frac{z^k}{\zeta^k}\right)&\mbox{otherwise}
\end{cases}.
$$
Note that $E(\zeta,e^{i\alpha}z)=E(e^{-i\alpha}\zeta,z)$ and that $\overline{E(\zeta,\bar z)}=E(\bar\zeta,z)$.
For a non-negative real number $t$ we denote by $[t]$ its integer part, that is, the integer $d$ such that
$d\leq t<d+1$.

\begin{theorem}[Hadamard Factorisation Theorem]
Let $f$ be an entire function of finite order $\rho$, let $\{z_n\}$ be its zeroes, with $0$ excluded
and all zeroes are repeated according to their multiplicity. Then, for every $\eps>0$,
\begin{equation}
\label{eq:had1}
\sum_{n\geq 1}\frac{1}{|z_n|^{\rho+\eps}}<+\infty
\end{equation}
and $f$ can be factored as
\begin{equation}
\label{eq:had2}
f(z)=z^me^{P(z)}\prod_{n\geq 1}E_{[\rho]}(z_n,z)
\end{equation}
where $m$ is an integer and $P$ is a polynomial of degree at most $\rho$.
\end{theorem}
Note that \eqref{eq:had1} implies that the product in \eqref{eq:had2} converges uniformly over compact sets.
Moreover, \eqref{eq:had2} shows that an entire function of finite order is essentially determined by its {\sl complex} zeroes.

Finally, note that in Hadamard's factorization theorem, one may replace $\rho$ by any $\rho'\geq\rho$.

\subsubsection{Reconstruction of a function from the modulus on two lines}

\begin{theorem}\label{th:2ent}
Let $\alpha_1,\alpha_2\in[0,2\pi[$ with $\alpha_1-\alpha_2\notin\pi\Q$.
Let $f$ and $g$ be two entire functions of finite order and assume that
\begin{equation}
\label{eq:2ent}
|g(xe^{i\alpha_1})|=|f(xe^{i\alpha_1})|\quad\mbox{and}\quad|g(xe^{i\alpha_2})|=|f(xe^{i\alpha_2})|
\quad\mbox{for all }x\in\R
\end{equation}
then $g=cf$ where $c$ is constant of modulus $1$.
\end{theorem}

\begin{proof}
If $f$ has a zero at $0$, then \eqref{eq:2ent} implies that $g$ also has a zero at $0$ of the same order. Without loss of generality, we may therefore assume that $f$ and $g$ have no zero at $0$.
For $k=1,2$, let us define the four entire functions of order $\rho$, $f_k(z)=f(ze^{i\alpha_k})$
and $g_k(z)=g(ze^{i\alpha_k})$ so that $|g_k(x)|=|f_k(x)|$ for all $x\in\R$.
This may be rewritten as
\begin{equation}
\label{eq.2theta}
g_k(z)\overline{g_k(\bar z)}=f_k(z)\overline{f_k(\bar z)}\quad\mbox{for all }z\in\R.
\end{equation}
But as this is an equality between two entire functions, it is valid for all $z\in\C$.

Now write the Hadamard factorisation of $f$ and $g$ as
$$
f(z)=e^{\sum_{j=0}^{[\rho]}a_j z^j}\prod_{z_l\in\mathcal{Z}_f\cap\mathcal{Z}_g}E_{[\rho]}(z_l,z)
\prod_{z_l\in\mathcal{Z}_f\setminus\mathcal{Z}_g}E_{[\rho]}(z_l,z)
$$
and
$$
g(z)=e^{\sum_{j=0}^{[\rho]}b_j z^j}\prod_{z_l\in\mathcal{Z}_f\cap\mathcal{Z}_g}E_{[\rho]}(z_l,z)
\prod_{z_l\in\mathcal{Z}_g\setminus\mathcal{Z}_f}E_{[\rho]}(z_l,z)
$$
where $\rho$ is greater than the orders of $f$ and $g$.
Notice that \eqref{eq.2theta} implies after simplification of the product over the common zeroes that
\begin{eqnarray*}
&&e^{2\sum_{j=0}^{[\rho]}\Re\bigl(a_je^{ij\alpha_k}\bigr)z^j}
\prod_{z_l\in\mathcal{Z}_f\setminus\mathcal{Z}_g}E_{[\rho]}(e^{-i\alpha_k}z_l,z)E_{[\rho]}(\overline{e^{-i\alpha_k}z_l},z)\\
&=&e^{2\sum_{j=0}^{[\rho]}\Re\bigl(b_je^{ij\alpha_k}\bigr)z^j}
\prod_{z_l\in\mathcal{Z}_g\setminus\mathcal{Z}_f}E_{[\rho]}(e^{-i\alpha_k}z_l,z)E_{[\rho]}(\overline{e^{-i\alpha_k}z_l},z).
\end{eqnarray*}
From this, we get the two following consequences:

--- first, $\dst\sum_{j=1}^{[\rho]}\Re\bigl(b_je^{ij\alpha_k}\bigr)z^j
=\sum_{j=1}^{[\rho]}\Re\bigl(a_je^{ij\alpha_k}\bigr)z^j$ for $k=1,2$ thus
$b_j=a_j$ for all $j\not=0$ and $\Re(b_0)=\Re(a_0)$ thus $b_0=a_0+i\gamma$, $\gamma\in\R$.

--- $z_l\in\mathcal{Z}_f\setminus\mathcal{Z}_g$ if and only if $e^{2i\alpha_k}\overline{z_l}\in \mathcal{Z}_g\setminus\mathcal{Z}_f$
for $k=1,2$ and vice versa. It results that
$(\mathcal{Z}_f\setminus\mathcal{Z}_g)\cup(\mathcal{Z}_g\setminus\mathcal{Z}_f)$ is stable under the reflections with respect to
the two lines $e^{-i\alpha_1}\R$ and $e^{-i\alpha_2}\R$ and thus under the rotation of angle $2(\alpha_1-\alpha_2)$.
But this angle is not a rational multiple of $\pi$ thus the orbit of a point is not a discrete set.
This is not compatible with the fact that
$\mathcal{Z}_f\setminus\mathcal{Z}_g\cup(\mathcal{Z}_g\setminus\mathcal{Z}_f)$ is a set of zeroes
of an entire function.
As a consequence $\mathcal{Z}_f\setminus\mathcal{Z}_g=\mathcal{Z}_g\setminus\mathcal{Z}_f=\emptyset$
that is $\mathcal{Z}_g=\mathcal{Z}_f$.

Summarising, we get $g(z)=e^{i\gamma}f(z)$ as announced.
\end{proof}

\begin{remark}\label{rq:2ent}
Note that it is enough to assume that $|g(xe^{i\alpha_j})|=|f(xe^{i\alpha_j})|$ for $x\in E_j$
where $E_j\subset \R$ is such that $E_je^{i\alpha_j}$ is a set of uniqueness for entire functions of order $\rho$. That is, $E_j$ is such that two entire functions of order $\rho$ that coincide on $E_je^{i\alpha_j}$ coincide everywhere. For instance, any set of positive measure would do.
\end{remark}

\begin{remark}\label{rq:2entb}
The key point in the previous proof is the Schwarz Reflection Principle and the fact that the (closed) orbit of a point under certain reflections can not be included in a zero-set of an entire function. One may therefore extend the previous theorem to several variable entire functions. We refrain from doing so as we have no specific application in mind.  
\end{remark}

\begin{remark}\label{rq:2entc}
In the case $\alpha_1-\alpha_2\in\pi\Q$, the previous proof still leads to a result if one further
restricts the properties of the zero set of the functions $f$ and $g$ {\it e.g.} to be in a strip. 

Indeed, the orbit of a point
in $\mathcal{Z}:=\mathcal{Z}_f\setminus\mathcal{Z}_g\cup(\mathcal{Z}_g\setminus\mathcal{Z}_f)$ is now a regular
polygon. Therefore, if one assumes that the zeroes of $f$ and $g$ are in a strip, then $\mathcal{Z}$ is bounded
therefore finite. It follows that $f$ and $g$ differ by at most a polynomial factor $f=Ph$, $g=Qh$
where $h$ is entire of finite order and $P,Q$ are polynomials. Moreover, the union of the zeros of $P$ and $Q$
is a finite union of regular polygons.
\end{remark}

\subsection{Some facts from Fourier analysis}

\subsubsection{Fourier transforms of compactly supported functions}
We will extensively make use of the classical Paley-Wiener Theorem:

\begin{theorem}[Paley-Wiener]
Let $f$ be an entire function such that $\vert f(z)\vert \leq K e^{2\pi\gamma |z|}$ for some $K \geq 0$ and $\gamma > 0$. If the restriction of $f$ to the real line is in $L^2(\mathbb{R})$ then there exists a function $u\in L^2(\R)$ with support included in $[-\gamma, \gamma]$ such that
$f=\ff u$.
\end{theorem}

Moreover, a theorem of Titchmarsch \cite{Ti} states that,

\begin{theorem}[Titchmarsch \cite{Ti}]
Let $u\in L^2(\R)$ be a compactly supported function so that $\ff u$ extends to an entire function
over $\C$. Let $z_n=r_ne^{i\theta_n}$, $0<r_1\leq r_2\leq\ldots$ be the nonzero zeroes of $\ff u$,
counted with multiplicity and arranged according to increasing modulus. Then
\begin{enumerate}
\renewcommand{\theenumi}{\roman{enumi}}
\item\label{cond:t2} $\dst\sum_{n=1}^{+\infty}\frac{1}{r_n}=+\infty$
but $\dst\sum_{n=1}^{+\infty}\frac{1}{r_n^{1+\eps}}<+\infty$ for every $\eps>0$,\\[2pt]

\item\label{cond:t3} $\dst\sum_{n=1}^{+\infty}\frac{|\sin\theta_n|}{r_n}<+\infty$
while $\dst\sum_{n=1}^{+\infty}\frac{\cos\theta_n}{r_n}$ converges (conditionally).
\end{enumerate}
\end{theorem}
Note that this theorem implies that $\ff u$ has infinitely many zeroes which can been seen
directly using Paley-Wiener's Theorem and Hadamard factorization. Indeed, an expression
like \eqref{eq:had2} can not be in $L^2(\R)$ if the product part is finite
since $P(z)e^{az}$ is not in $L^2(\R)$ when $P$ is a non-zero polynomial.

\subsubsection{The Shannon-Whittacker sampling theorem}
Results in this section may be found in \cite{BZ} and \cite{Un}
and in the numerous references therein. They may also
be found in numerous books on Fourier analysis, although the oversampling 
formula is sometimes slightly erroneous.

Let us recall that the $\sinc$ function is defined by
$$
\sinc x=\begin{cases}\frac{\sin x}{x}&\mbox{if }x\not=0\\
1&\mbox{if }x=0\end{cases}.
$$
On $\R^d$ one defines $\sinc(x)=\sinc(x_1)\cdots\sinc(x_d)$. 
%
%
%

\begin{theorem}[Shannon-Whittacker Sampling Formula (with oversampling]\label{th:oversample}
Let $\sigma$ be a positive real number and let $\dst h\leq\frac{1}{\sigma}$.
Let $\gamma\in\cc^j(\R^d)$ be a non-negative even function such that
$\gamma(\xi)=0$ for $|\xi|\geq 1$ and normalised by
$$
\int_{\R^d}\gamma(\xi)\,\mathrm{d}\xi=1.
$$
Let $\ffi\in L^2(\R^d)$ be such that its Fourier transform be supported in $[-\sigma,\sigma]^d$.
Then $\ffi$ is uniquely defined from its samples $\ffi(hk/2)$, $k\in\Z^d$ by the formula
\begin{equation}
\label{eq:shannon4}
\ffi(x)=(2\tau_+)^d\sum_{k\in\Z^d}\ffi\left(\frac{hk}{2}\right)
\widehat{\gamma}\left(\tau_-\left(x-\frac{hk}{2}\right)\right)
\sinc 2\pi\tau_+\left(x-\frac{hk}{2}\right)
\end{equation}
where $\dst\tau_\pm=\frac{1\pm\sigma h}{2h}$.
\end{theorem}

This result is usually stated without oversampling, that is the case $\sigma h=1$ in which case it reduces to the usual Shannon-Whittacker Formula:
\begin{equation}
\ffi(x)=(2\sigma)^d\sum_{k\in\Z^d}\ffi\left(\frac{hk}{2}\right)\sinc\frac{2\pi}{h}\left(x-\frac{hk}{2}\right).
\label{eq:shannon1}
\end{equation}
The use of oversampling improves the convergence properties of the series
\eqref{eq:shannon1}. Note also that various error bounds are known and can be found
in text books like \cite{Hi,Za}: aliasing ($u$ is only approximatively band-limited) jittering (the samples
are only taken approximatively at $hk$), truncation (the series \eqref{eq:shannon1} is truncated)...

\subsection{Preliminaries on the fractional Fourier transform}

\subsubsection{Definitions}\label{sec:deffrac}
For $\alpha\in\R\setminus\pi\Z$, let $c_\alpha=\dst\frac{\exp\frac{i}{2}\left(\alpha-\frac{\pi}{2}\right)}{\sqrt{|\sin\alpha|}}$
be a square root of $1-i\cot\alpha$. For $u\in L^1(\R^d)$ and $\alpha\notin\pi\Z$, define
\begin{eqnarray}
\ff_\alpha u(\xi)&=&c_\alpha^d e^{-i\pi|\xi|^2\cot\alpha}\int_{\R^d}
u(t)e^{-i\pi |t|^2\cot\alpha}e^{-2i\pi \scal{t,\xi/\sin\alpha}}\mbox{d}t\notag\\
&=&c_\alpha^d e^{-i\pi|\xi|^2\cot\alpha}\ff[u(t)e^{-i\pi |t|^2\cot\alpha}](\xi/\sin\alpha)
\label{eq:deffrfo}
\end{eqnarray}
while for $k\in\Z$, $\ff_{2k\pi} u=u$ and $\ff_{(2k+1)\pi}u(\xi)=u(-\xi)$.
This transformation has the following properties~:
\begin{enumerate}
\item\label{prop:fa1} $\dst\int_{\R^d}\ff_\alpha u(\xi)\overline{\ff_\alpha v(\xi)}\mbox{d}\xi=\int_{\R^d} u(t)\overline{v(t)}\mbox{d}t$
which allows to extend $\ff_\alpha$ from $L^1(\R^d)\cap L^2(\R^d)$ to
$L^2(\R^d)$ as a unitary operator on $L^2(\R^d)$;

\item\label{prop:fa2} $\ff_\alpha\ff_\beta=\ff_{\alpha+\beta}$;

\item\label{prop:fa3} if $u_{a,\omega}(t)=e^{-2i\pi\scal{\omega,t}}u(t-a)$, $a,\omega\in\R^d$, then
$$
\ff_\alpha u_{a,0}(\xi)=\ff_\alpha u(\xi+a\cos\alpha)e^{-i\pi |a|^2\cos\alpha\sin\alpha -2i\pi \scal{a,\xi}\sin\alpha};
$$
and
$$
\ff_\alpha u_{0,\omega}(\xi)=\ff_\alpha u(\xi+\omega\sin\alpha)e^{i\pi|\omega|^2\cos\alpha\sin\alpha+2i\pi\scal{\omega,\xi}\cos\alpha}.
$$
\end{enumerate}

Further, let us define the Hermite basis functions on $\R$ by
$$
h_k(t)=\frac{2^{1/4}}{\sqrt{k!}}\left(-\frac{1}{\sqrt{2\pi}}\right)^ke^{\pi t^2} \left(\frac{\mbox{d}}{\mbox{d}t}\right)^ke^{-2\pi t^2},\quad k\in\N.
$$
Hermite functions are then defined on $\R^d$ by tensorisation
$$
h_{k_1,\ldots,k_d}(x_1,\ldots,x_d)=h_{k_1}(x_1)\cdots h_{k_d}(x_d).
$$
It is well known that $(h_k)_{k\in\N^d}$ is an orthonormal basis of $L^2(\R^d)$
and that $h_k$ can be expressed as $h_k(x)=H_k(x)e^{-\pi |x|^2}$ with $H_k$ a polynomial of degree $|k|=k_1+\cdots+k_d$.

The fractional Fourier transform may alternatively be defined as
\begin{equation}
\label{eq:hermite}
\ff_\alpha[u]=\sum_{(k_1,\ldots,k_d)\in\N^d}e^{-i\alpha(k_1+\cdots+k_d)}
\scal{u,h_{k_1,\ldots,k_n}}h_{k_1,\ldots,k_n}.
\end{equation}

\begin{remark}\label{remfatriv}
Note that the fractional Fourier transform $\ff_\alpha[u]$
($\alpha\notin\pi\Z$) is obtained from the ordinary Fourier transform in the following way:
take $u$, multiply it by a chirp $e^{-i\pi |t|^2\cot\alpha}$, take the ordinary Fourier transform,
re-scale and multiply by a chirp $e^{i\pi |t|^2\cot\alpha}$.
As a consequence, it is often straightforward to translate results on the ordinary Fourier transform
({\it e.g.} Paley-Wiener and Titshmarsh's Theorems that will be used below) into results for fractional Fourier transforms.
\end{remark}

\subsubsection{The Phase Retrieval Problem for the fractional Fourier transform
for compactly supported functions}\label{sec:phafr}
In this section, for the sake of completeness, we will solve the following Phase Retrieval Problem:

\medskip

\noindent{\bf Problem. (FrFT Phase Retrieval Problem)}\\
{\sl Let $\alpha\in\R\setminus\pi\Z$. Given $u\in L^2(\R)$ with compact support,
find all $v\in L^2(\R)$ with compact support such that,
for every $\xi\in\R$, $|\ff_\alpha v(\xi)|=|\ff_\alpha u(\xi)|$.}

\medskip

To start with, let us write $u_\alpha(t)=e^{-i\pi t^2\cot \alpha}u(t)$
so that $u_\alpha\in L^2(\R)$ with compact support and
$\ff_\alpha u(\xi)=c_\alpha e^{-i\pi\xi^2\cot\alpha}\ff u_\alpha(\xi/\sin\alpha)$. It follows that
$|\ff_\alpha v(\xi)|=|\ff_\alpha u(\xi)|$ is equivalent to
\begin{equation}
\label{eq:hol}
|\ff v_\alpha(\xi)|^2=|\ff u_\alpha(\xi)|^2.
\end{equation}

Now Paley-Wiener's theorem implies that $\ff u_\alpha$ and $\ff v_\alpha$
are entire functions of exponential type.
There are two consequences.
\smallskip

--- Rewriting \eqref{eq:hol} as
\begin{equation}
\label{eq:hol2}
\ff v_\alpha(\xi)\overline{\ff v_\alpha(\bar\xi)}=\ff u_\alpha(\xi)\overline{\ff u_\alpha(\bar\xi)}, 
\end{equation}
we have an equality of entire functions, so that its validity extends from $\xi\in\R$ to $\xi\in\C$.

\smallskip

--- According to Hadamard's Factorisation Theorem, $\ff u_\alpha$ may be written in the form
\begin{equation}
\label{eq:had}
\ff u_\alpha(\xi)=\xi^{k_u} e^{a_u\xi+b_u}\prod_{\zeta\in\zzz_{u}}E_1(\zeta,\xi),
\end{equation}
where $k_u$ is an integer, $a_u,b_u\in\C$ and $\zzz(u):=\{r_ne^{i\theta_n},n\in\N\}$ the set of complex zeroes of
$\ff_\alpha(u)$ (all depending on $\alpha$).
Moreover, Titchmarsch's theorem implies that the product in \eqref{eq:had} converges uniformly over compact sets. A similar representation holds for $\ff v_\alpha$. 

\smallskip

It follows then from
\eqref{eq:hol2} that $k_v=k_u$,
$\mbox{Re}\,a_v=\mbox{Re}\,a_u$, $\mbox{Re}\,b_v=\mbox{Re}\,b_u$
and $\zeta\in \zzz_v$ if and only if either $\zeta\in\zzz_u$ or
$\bar\zeta\in\zzz_u$. In other words, there exist $\omega\in\R$,
$c\in\C$ with $|c|=1$ and 
for each $n\in\N$, a choice $\eps_n=\pm1$ such that
$$
\ff v_\alpha(\xi)=ce^{i\omega\xi} \xi^{k_u}e^{a_u\xi+b_u}
\prod_{n=1}^{+\infty}E_1(r_ne^{i\eps_n\theta_n},\xi),
$$
which has to be compared to \eqref{eq:had}.
The convergence of this infinite product to an entire function of order $1$
is guaranteed by Titchmarch's Theorem while the fact that one recovers
an $L^2$ function is guaranteed by the fact that $|\ff v_\alpha|=|\ff u_\alpha|$
on $\R$. One then recovers $v$ by taking the inverse Fourier transform and multiplying the result by $e^{i\xi^2\cot\alpha}$.

\begin{remark}
The above result is directly inspired by the result on the usual phase retrieval theorem ($\alpha=\pi/2$) for which the corresponding result was first proved by Walter in \cite{Wa}, though in a less precise mathematical form. For proper mathematical justification of Walter's Theorem, we refer {\it e.g.} to \cite{Hu}.

The choice of $\eps_n=-1$ amounts to changing the complex root $r_ne^{i\theta_n}$
of $\ff_\alpha u$ into its complex conjugate and is called \emph{zero-flipping}
in the engineering literature.
\end{remark}

\begin{remark}
The same proof applies as soon as we impose conditions on $u$ and $v$ that force their fractional Fourier
transforms $\ff_\alpha[u]$ and $\ff_\alpha[v]$ to be (extended to) entire functions of finite type.
For instance, we could impose that $u$ and $v$ are bounded by gaussians $|u(x)|,|v(x)|\leq Ce^{-\lambda|x|^2}$
for some $\lambda>0$.
\end{remark}

\subsection{Preliminaries on the ambiguity function}\label{sec:amb}
Let us recall that the ambiguity function of $u,v\in L^2(\R)$ is defined by
$$
A(u,v)(x,y)=\int_{\R^d} u\left(t+\frac{x}{2}\right)\overline{v\left(t-\frac{x}{2}\right)}
e^{-2i\pi \scal{t,y}}\mbox{d}t.
$$
We will simply denote $A(u)=A(u,u)$. Note that $A(u)(0,y)=\ff[|u|^2](y)$.

The reader may be more acquainted wit the following closely related transforms:
the short-time Fourier transform, also known as the
windowed Fourier transform, defined by $S_vu(x,y)=e^{-i\pi\scal{x,y}}A(u,v)(x,y)$ and the
Wigner transform
$$
W(u,v)(\eta,\xi)=\int_{\R^d} u\left(\xi+\frac{t}{2}\right)\overline{v\left(\xi-\frac{t}{2}\right)}
e^{2i\pi \scal{t,\eta}}\mbox{d}t
$$
which is the inverse Fourier transform of $A(u,v)$ in $\R^{2}$. Further
$$
W(u,v)(\eta,\xi)=2^dA(u,\check v)(2\xi,2\eta)
$$
where $\check v(t)=v(-t)$.

Let us now list the properties that we need. 
They are all well-known {\it e.g.} \cite{Al,AT,Wi}:
\begin{enumerate}
\item\label{prop:amb1} $A(u,v)\in L^2(\R^{2d})$ with $\norm{A(u,v)}_{L^2(\R^{2d})}=\norm{u}_2\norm{v}_2$ and is continuous, going to $0$ at infinity;

\item\label{prop:amb2}  $\norm{A(u)}_\infty=A(u)(0,0)=\norm{u}_2^2$;

\item\label{prop:amb2a} $A(v)=A(u)$ if and only if there
exists $c\in\C$ with $|c|=1$ such that $v=cu$

\item\label{prop:amb3b} time-frequency shifts are transformed into time-frequency shifts:
for $a,\omega\in\R^d$ let $u_{a,\omega}(t)=e^{2i\pi\scal{\omega,t}}u(t-a)$, then
$$
A(u_{a,\omega},v_{b,\eta})(x,y)=
e^{i\pi\bigl(\scal{\omega+\eta,x}+\scal{a+b,y-\omega+\eta}\bigr)}A(u,v)(x-a+b,y-\omega+\eta);
$$

\item\label{prop:amb4} The fractional Fourier transform rotates the variables:
$$
A(\ff_\alpha u,\ff_\alpha v)(x,y)=A(u,v)(x\cos\alpha-y\sin\alpha,x\sin\alpha+y\cos\alpha).
$$
In particular,
\begin{equation}
\label{eq:key}
A(u)(-y\sin\alpha,y\cos\alpha)=A(\ff_\alpha u)(0,y)=\ff[|\ff_\alpha u|^2](y).
\end{equation}

\end{enumerate}

Property \ref{prop:amb4} was first proved in \cite{Wi} but the use of the fractional Fourier transform 
may have been unnoticed as it is defined via Hermite functions.
With the above expression of the fractional Fourier transform, this result seems to have first appeared in
\cite{Na} and has been rediscovered several times {\it e.g.} \cite{Al,Lo}.

This property also justifies the interpretation of the fractional Fourier transform as a rotation of a function in the ``time-frequency plane''.

\section{Reconstruction from phase-less windowed Fourier transform or wavelet transform with well chosen windows}\label{sec:window}

Before going on with the main scope of this paper, let us illustrate how Theorem \ref{th:2ent}
can be applied to the reconstruction of signals from phase-less measurements of particular
windowed transforms.

\subsection{Reconstruction from a windowed Fourier transform}

\begin{proposition}\label{prop:window}
Let $\gamma(t)=e^{-\pi t^2}$ and let $\alpha\in\R\setminus\pi\Q$.
Let $u,v\in L^2(\R)$, and assume that
\begin{equation}
\label{eq:window}
|v*\gamma (t)|=|u*\gamma (t)|\quad\mathrm{and}\quad
|\ff_\alpha[v]*\gamma (t)|=|\ff_\alpha[u]*\gamma (t)|
\end{equation}
for $t\in \R$. Then there exists $c\in\C$ with $|c|=1$ such that $v=cu$.
\end{proposition}

\begin{proof}
First note that the windowed Fourier transform with window $\gamma$ is given by
\begin{eqnarray*}
S_\gamma u(x,y)&=&\int_\R u(t)\gamma(t-x)e^{2i\pi ty}\,\mbox{d}t\\
&=&e^{-\pi y^2+2i\pi xy}\int_\R u(t)e^{-\pi(t-z)^2}\,\mbox{d}t
=e^{-\pi y^2+2i\pi xy}F_u(z)
\end{eqnarray*}
where $z=x+iy$ and $\dst F_u(z)=\int_\R u(t)e^{-\pi(t-z)^2}\,\mbox{d}t$.

Further $u*\gamma(t)=S_\gamma u(t,0)$ and, using the fact that the fractional Fourier transform rotates the
variables of the ambiguity function (Property \ref{prop:amb4} of the ambiguity function)
and leaves the Gaussian invariant $\ff_\theta\gamma=\gamma$, we get that 
\begin{eqnarray*}
S_\gamma u(t\cos\alpha,t\sin\alpha)&=&e^{-i\pi t^2\cos\alpha\sin\alpha}A(u,\gamma)(t\cos\alpha,t\sin\alpha)\\
&=&e^{-i\pi t^2\cos\alpha\sin\alpha}A(\ff_\alpha u,\ff_\alpha\gamma)(t,0)\\
&=&e^{-i\pi t^2\cos\alpha\sin\alpha}A(\ff_\alpha u,\gamma)(t,0)
\end{eqnarray*}
since $\ff_\alpha\gamma=\gamma$. Unwinding the computation, we obtain
$$
S_\gamma u(t\cos\alpha,t\sin\alpha)=e^{-i\pi t^2\cos\alpha\sin\alpha}
S_\gamma\ff_\alpha u(t,0)=e^{-i\pi t^2\cos\alpha\sin\alpha}\ff_\alpha[u]*\gamma(t).
$$

From the above computations, we see that \eqref{eq:window} implies that
$|F_u(t)|=|F_v(t)|$ and $|F_u(te^{i\alpha})|=|F_v(te^{i\alpha})|$ for $t\in\R$.
It remains to notice that, for $u;v\in L^2(\R)$, $F_u,F_v$ are entire functions
of order $2$ so that Theorem \ref{th:2ent} gives the result.
\end{proof}

\begin{remark}
Note that it is enough to have \eqref{eq:window} on a set of positive measure.
\end{remark}

\begin{remark}
By using the Paley-Wiener Theorem for Schwatz class distributions, one easily sees that the previous
result stays true if $u,v\in\ss'(\R)$, provided one extends the definition of the fractional Fourier transform
to $\ss'$ as well.
\end{remark}

\subsection{Reconstruction from a wavelet transform}
Let us recall that the wavelet transform with window $\psi\in L^2(\R)$ is defined for $u\in L^2(\R)$ by
$$
\mathcal{W}_\psi(u)(a,\tau)=\sqrt{a}\int_\R u(t)\overline{\psi\bigl(a(t-\tau)\bigr)}\,\mbox{d}t.
$$
A simple computation using Parseval's identity shows that
$$
\mathcal{W}_\psi(u)(a,\tau)=\widetilde{\mathcal{W}}_{\hat\psi}(\hat u)(1/a,\tau/a)
$$
where
$$
\widetilde{\mathcal{W}}_{\hat\psi}(\hat u)(b,\xi)=\sqrt{b}\int_\R \hat u(t)\overline{\hat\psi(bt)}e^{2i\pi t\xi}\,\mbox{d}t.
$$

Now, let $\rho,\beta>0$, and let us define 
$\psi(x)=\psi_\pm(x)=\dst\frac{1}{2\pi(\rho\mp ix)}$ so that
$\hat\psi_+(x)=\chi_{[0,+\infty)}e^{-2\pi \rho t}$ and 
$\hat\psi_-(x)=\chi_{(-\infty,0]}e^{2\pi \rho t}$.
If $\hat u$ is supported in $(-\infty,\beta]$
--- resp. in $[-\beta,+\infty)$--- then
\begin{equation}
\label{eq:fdet}
F_{u,+}(\xi):=\widetilde{\mathcal{W}}_{\hat\psi_+}(\hat u)(1,\xi)=\int_0^{\beta}\hat u(t)e^{-2\pi\rho t}e^{2i\pi t\xi}\,\mbox{d}t
\end{equation}
--- resp. $F_{u,-}(\xi):=\widetilde{\mathcal{W}}_{\hat\psi_+}(\hat u)(1,\xi)$---
is an entire function of order $1$. But, if $\alpha\in[0,\pi/2)$ and if $\xi\in\R$ such that
$1+\rho^{-1}\xi\sin\alpha>0$,
\begin{eqnarray*}
F_{u,+}(\xi e^{i\alpha})&=&\int_0^{\beta}\hat u(t)e^{-2\pi\rho(1+\rho^{-1}\xi\sin\alpha) t}
e^{2i\pi t\xi\cos\alpha}\,\mbox{d}t\\
&=&(1+\rho^{-1}\xi\sin\alpha)^{-1/2}\widetilde{\mathcal{W}}_{\hat\psi}(\hat u)(1+\rho^{-1}\xi\sin\alpha,\xi\cos\alpha).
\end{eqnarray*}

Now let $v$ be such that $\hat v$ is also supported in $(-\infty,\beta]$
and assume that  $|\mathcal{W}_{\psi_+}(v)(1,\tau)|=|\mathcal{W}_{\psi_+}(u)(1,\tau)|$
and, if $1+\rho^{-1}\xi\sin\alpha>0$,
$$
\abs{\mathcal{W}_{\psi_+}(v)\left(\frac{\rho}{\rho+\xi\sin\alpha},\frac{\rho\xi\cos\alpha}{\rho+\xi\sin\alpha}\right)}
=\abs{\mathcal{W}_{\psi_+}(u)\left(\frac{\rho}{\rho+\xi\sin\alpha},\frac{\rho\xi\cos\alpha}{\rho+\xi\sin\alpha}\right)}.
$$
Expressing this in terms of $\widetilde{\mathcal{W}}$ and then $F$ as above, one easily checks that
$|F_{v,+}(\xi)|=|F_{u,+}(\xi)|$ for $\xi\in\R$
and $|F_{v,+}(\xi e^{i\alpha})|=|F_{u,+}(\xi e^{i\alpha})|$ for $\xi$ such that $1+\rho^{-1}\xi\sin\alpha>0$.
It follows from Theorem \ref{th:2ent} and Remark \ref{rq:2ent} that, 
$F_{v,+}=c_+F_{u,+}$ where $|c_+|=1$.
But then, taking Fourier transforms, we obtain
$$
\hat u(t)\chi_{[0,+\infty)}(t)=e^{2\pi \rho t}\ff[F_{u,+}](t)
\quad\mbox{and}\quad
\hat v(t)\chi_{[0,+\infty)}(t)=e^{2\pi \rho t}\ff[F_{v,+}](t)
$$ 
thus $\hat v\chi_{[0,+\infty)}=c_+\hat u\chi_{[0,+\infty)}$.

Replacing $\psi_+$ and $\psi_-$ and assuming that $\hat u,\hat v$ are supported in
$[-\beta,+\infty)$ we obtain $\hat v\chi_{(-\infty,0]}=c_-\hat u\chi_{(-\infty,0]}$
with $|c_-|=1$.

Summarising (and simplifying the condition), we obtain the following:

\begin{theorem}
Let $\rho>0$ and let $\psi_\pm$ be defined on $\R$ as $\psi_\pm(x)=\dst\frac{1}{2\pi(\rho\mp ix)}$.
Let $u,v\in L^2(\R)$ have compactly supported Fourier transforms ({\it i.e.} be band-limited)
and such that $|\mathcal{W}_{\psi_\pm}(v)|=|\mathcal{W}_{\psi_\pm}(u)|$.
Then there exist $c_\pm\in\C$ with $|c_\pm|=1$ such that
$$
\hat v(\xi)=\begin{cases}c_+\hat u(\xi)\mbox{ on }\R^+\\
c_-\hat u(\xi)\mbox{ on }\R^-
\end{cases}.
$$
\end{theorem}

\begin{remark}
Note that if we further assume that $\hat u,\hat v$ are continuous, then $c_+=c_-$.
\end{remark}

\section{Reconstruction of a function from multiple fractional Fourier transform intensities}
\label{sec:reconstruction}

\subsection{Reconstruction from the modulus of many fractional Fourier transforms}\ \\
We will now prove the following result:

\begin{theorem}\label{mainthm}
In the following cases, exact reconstruction can be obtained.
\begin{enumerate}
\item Let $u,v\in L^2(\R)$ such that, for every $\alpha\in[-\pi/2,\pi/2]$, $|\ff_\alpha v|=|\ff_\alpha u|$. Then there exists $c\in\C$ with $|c|=1$ such that $v=cu$.

\item\label{th-point2}  Let $\tau\subset[-\pi/2,\pi/2]$ be either a set of positive measure or a set with an accumulation point
$\alpha_0\not=0$.
Let $u,v\in L^2(\R)$ with compact support such that,
for every $\alpha\in\tau$, $|\ff_\alpha v|=|\ff_\alpha u|$.
Then there exists $c\in\C$ with $|c|=1$ such that $v=cu$.

\item Let $a>0$ and define $(\alpha_k)_{k\in\Z}$
by $\alpha_0=\pi/2$ and, for $k\in\Z\setminus\{0\}$, $\alpha_k=\dst\arctan\frac{a^2}{k}$. Then,
given $u,v\in L^2(\R)$ with compact support included in $[-a,a]$,
if $|\ff_{\alpha_k} v|=|\ff_{\alpha_k} u|$ for every $k\in\Z$,
then there exists $c\in\C$ with $|c|=1$ such that $v=cu$.
\end{enumerate}
\end{theorem}

\begin{proof}
From \eqref{eq:key}, if $|\ff_\alpha v|=|\ff_\alpha u|$ for every $\alpha\in[-\pi/2,\pi/2]$,
then, from \eqref{eq:key}, we get $A(v)=A(u)$.
Property \ref{prop:amb2a} of the ambiguity function then implies that
$v=cu$ for some complex number with $|c|=1$.

\begin{remark}
As $A(u)$ is continuous, it is enough to have $\alpha$ in a dense subset of $[-\pi/2,\pi/2]$
for the previous result to hold.
\end{remark}

For the second part of the theorem, \eqref{eq:key} implies that 
$$
A(v)(-y\sin\alpha,y\cos\alpha)=A(u)(-y\sin\alpha,y\cos\alpha)
$$
for every $y$ and every $\alpha\in\tau$. In particular, if we fix $x$ and denote by
$$
\tau_x:=\{-x/\tan\alpha,\ \alpha\in \tau\},
$$
then $A(v)(x,\xi)=A(u)(x,\xi)$ for every $\xi\in \tau_x$. 
But, if $\tau$ is of positive measure (resp. has an accumulation point $\alpha_0\not=0$), then so is $\tau_x$
(resp. has an accumulation point at $-x/\tan\alpha_0$ --- at $0$ if $\alpha_0=\pm\pi/2$) .
On the other hand $A(u)(x,\xi)$ is the Fourier transform of
$\ffi_x$ defined by $\ffi_x(t)=u\left(t+\frac{x}{2}\right)\overline{u\left(t-\frac{x}{2}\right)}$.
This function is of compact support so that $A(u)(x,\xi)$ is an entire function
in the $\xi$ variable. The same is true for $A(v)(x,\xi)$.
Finally, if two entire functions agree on a set of positive measure
(resp. with an accumulation point), they agree everywhere.

\medskip

For the last part of the theorem, note that
$\ffi_x$ is an $L^1$ function supported in the interval $\dst\ent{-a+\frac{|x|}{2},a-\frac{|x|}{2}}$ when $|x|<2a$ and is $0$ for $|x|\geq 2a$.

But then, from the Shannon-Whittaker Formula, the Fourier transform $A(u)(x,y)$ of $\ffi_x$
may be reconstructed from its samples. More precisely,
\begin{eqnarray}
A(u)(x,y)&=&\widehat{\ffi_x}(y)
=\sum_{k\in\Z}\widehat{\ffi_x}\left(\frac{h_xk}{2}\right)\,\mbox{sinc}
\frac{2\pi}{h_x}\left(y+\frac{h_xk}{2}\right)\nonumber\\
&=&\sum_{k\in\Z}A(u)\left(x,\frac{h_xk}{2}\right)\,\mbox{sinc}\frac{2\pi}{h_x}\left(y+\frac{h_xk}{2}\right)
\label{eq:shanau}
\end{eqnarray}
provided $|h_x|\leq\dst\frac{1}{a-|x|/2}$. A similar expression holds for $A(v)(x,y)$.

Now recall  from \eqref{eq:key} that
$$
\ff[|\ff_{-\alpha} u|^2](\xi)=A(\ff_{-\alpha} u)(0,\xi)=A(u)(-\xi\sin\alpha,\xi\cos\alpha).
$$
In particular, if we chose $h_x$ of the form $h_x=\gamma x$ then,
$$
A(u)\left(x,\frac{h_xk}{2}\right)=A(u)\left(x,\frac{\gamma k}{2}x\right)
=\ff[|\ff_{-\alpha_k} u|^2](-x/\sin\alpha_k)
$$
where $\alpha_k=\dst\arctan\frac{2}{k\gamma}$ ($\alpha_0=-\pi/2$). 

As we assumed that $|\ff_{-\alpha_k} v|=|\ff_{-\alpha_k} u|$ for every $k$,
it follows that $A(v)\left(x,\frac{h_xk}{2}\right)=A(u)\left(x,\frac{h_xk}{2}\right)$
for every $x$ and $k$. From \eqref{eq:shanau} we then deduce that $A(v)=A(u)$ everywhere thus
$v=cu$ with $c\in\C$, $|c|=1$.

It remains to choose $\gamma$ so that $\dst\abs{\gamma x}\leq\frac{1}{a-|x|/2}$ for $0\leq x\leq 2a$.
For this, it is enough to find the point on the hyperbola $y=1/(a-x/2)$ at which the tangent goes 
through the origin. Easy calculus shows that this point is $x=a$, $y=2/a$, thus
any $\gamma\leq\dst\frac{2}{a^2}$ will do ({\it see} Figure \ref{fig:nyquist}).
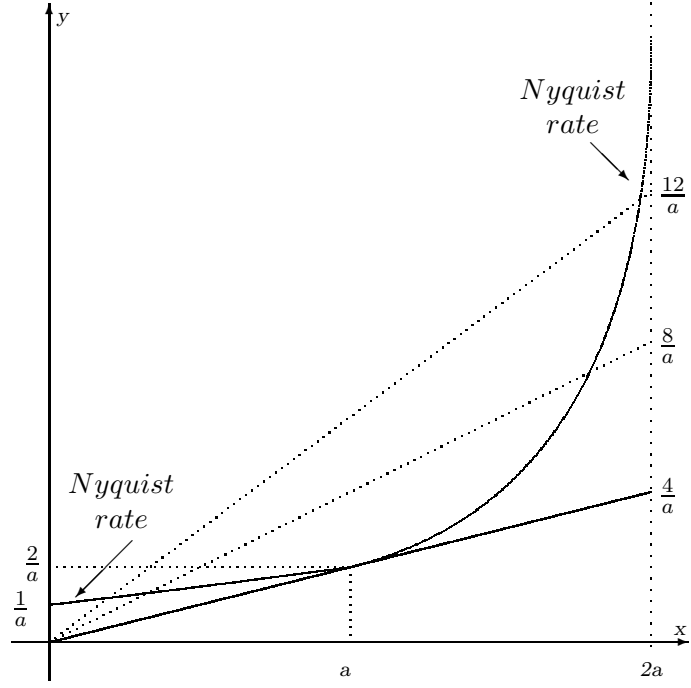
\begin{figure}
\begin{center}
\setlength{\unitlength}{0.5cm}
\begin{picture}(18,18)
\put(1,0){\vector(0,1){18}}
\put(0,1){\vector(1,0){18}}
\dottedline{0.3}(17,0.9)(17,18)
\dottedline{0.01}(1,1)(17,5)
\put(17.2,4.7){$\frac{4}{a}$}
\put(17.2,8.7){$\frac{8}{a}$}
\dottedline{0.2}(1,1)(17,9)
\put(17.2,12.7){$\frac{12}{a}$}
\dottedline{0.2}(1,1)(17,13)
\qbezier[200](1,2)(3.7,2.3)(9,3)
\qbezier[300](9,3)(17,5)(17,17)
\put(0,1.7){$\frac{1}{a}$}
\put(0.3,2.9){$\frac{2}{a}$}
\dottedline{0.2}(1,3)(9,3)(9,1)
\put(8.7,0.1){\sl\scriptsize a}
\put(16.7,0.1){\sl\scriptsize 2a}
\put(17.6,1.2){\sl\scriptsize x}
\put(1.2,17.5){\sl\scriptsize y}
\put(13.5,15){$\begin{matrix}Nyquist\\ rate\end{matrix}$}
\put(15.4,14.3){\vector(1,-1){1}}
\put(1.5,4.5){$\begin{matrix}Nyquist\\ rate\end{matrix}$}
\put(3.2,3.7){\vector(-1,-1){1.5}}
\end{picture}
\caption{The Nyquist rate and the sampling lines for $A(u)$}
\label{fig:nyquist}
\end{center}
\end{figure}

%
\end{proof}

\begin{remark}
The second part can be extended in many ways. For instance, it is enough
to assume that $u$ and $v$ have fast enough decrease to force
$A(u)$ and $A(v)$ to be entire functions in the second variable. For instance, this is the case
if $u$ and $v$ have a Gaussian bound $|u(x)|\leq Ce^{-\lambda |x|^2}$, $\lambda>0$.
\end{remark}

The last part of the theorem can be improved if one increases the number of measures in the following way:

\medskip

Let $b>a>0$ and define $(\alpha_k)_{k\in\Z}$
by $\alpha_0=\pi/2$ and, for $k\in\Z\setminus\{0\}$, $\alpha_k=\arctan\frac{k}{b^2}$. Assume that,
$u\in L^2(\R)$ has compact support included in $[-a,a]$, and
assume that we know
$|\ff_{-\alpha_k} u|$ for every $k\in\Z$. 
Then $u$ can be reconstructed in the following way:

Let $\gamma$ be a function that  satisfies the conditions of Theorem 
\ref{th:oversample}, that is, $\gamma$ is a smooth even 
non-negative function supported in $[-1,1]$ such that $\int\gamma=1$ and let 
$\hat\gamma$ be its Fourier Transform. Let 
$$
\psi_k(x)=\dst\ff[|\ff_{-\alpha_k} u|^2]
\left(-x\sin\alpha_k+\frac{|x|^2k\cos\alpha_k}{b^2}\right)
$$
and
$h_x=2x/b^2$.
Then, 
$$
A(u)(x,y)=\frac{b+2x}{2b}\sum_{k\in\Z}\psi_k(x)
\hat\gamma\left(\frac{b^2-2xb}{4x}\left(y-\frac{xk}{b^2}\right)\right)\mbox{sinc}
\,\frac{b^2+2xb}{4x}\left(y-\frac{xk}{b^2}\right).
$$
It is then enough to invert $A(u)$.

There are several issues in the practical implementation of this formula. We need to estimate
the error when reconstructing $A(u)$ (aliasing, jittering and truncation) and the error
when reconstructing $u$ from an approximation of $A(u)$.
We postpone this study to a future paper.

\subsection{Reconstruction of finite combinations of Hermite functions}\ \\
In this section, we will restrict our attention to signals $u\in L^2(\R^d)$
of the form:

\begin{definition}
Let $u\in L^2(\R^d)$ and $N$ a non-negative integer. We will say that
$u$ is a Hermite function of degree $N$, if $u$ is of the form
$u(t)=P(t)e^{-\pi|t|^2}$ where $P$ is a polynomial of degree $N$.
\end{definition}

Recall that we defined the elements of the Hermite basis $h_{k_1,\ldots,k_d}$ in Section \ref{sec:deffrac}. A Hermite function of degree $N$ may thus be expressed as
$$
u(t)=\sum_{\begin{matrix}
k_1,\ldots,k_d\in\N\\
k_1+\cdots+k_d\leq N
\end{matrix}}c_{k_1,\ldots,k_d}h_{k_1,\ldots,k_d}(t)
$$
with $c_{k_1,\ldots,k_d}\in\C$. Further, recall that the $h_{k_1,\ldots,k_d}$'s form an orthonormal basis of $L^2(\R^d)$
so that every $L^2$-function is well approximated by a Hermite function of degree $N$,
provided $N$ is big enough.

We may now prove the following:

\begin{theorem}\label{th:hermite}
Let $d\geq1$, $M,N\geq 0$ be integers and let $\alpha\not=\beta\in\R$ and $\gamma=\beta-\alpha$.
Assume that either $\gamma\notin\Q\pi$ or $\gamma=\dst\frac{p}{q}\pi$ with $p,q$ mutually prime integers and $q>\min(M,N)d$.

Let $u,v\in L^2(\R^d)$ be Hermite functions of degree $N$
and $M$ respectively. Assume that $|\ff_\alpha v|=|\ff_\alpha u|$ and that 
$|\ff_\beta v|=|\ff_\beta u|$, then there exists $c\in\C$ with $|c|=1$ such that $v=cu$.
\end{theorem}

The condition on $\gamma$ can also be reformulated as follows:
none of $e^{ij\gamma}$, $j=1$,..., $\min(M,N)d$ is real.

\begin{proof} Without loss of generality, we may assume that $\min(M,N)=N$.
For simplicity, let us give the proof only in the one-dimensional case.
Write
$$
u(t)=e^{-\pi t^2}\sum_{k=0}^Nc_kH_k(t)\quad,\quad
v(t)=e^{-\pi t^2}\sum_{k=0}^Md_kH_k(t)
$$
with $c_N\not=0$ and $d_M\not=0$.
It follows that, for $\theta\in\R$,
$$
\ff_\theta[u](t)=e^{-\pi t^2}\sum_{k=0}^Nc_ke^{-ik\theta}H_k(t)\quad,\quad
\ff_\theta[v](t)=e^{-\pi t^2}\sum_{k=0}^Ne^{-ik\theta}d_kH_k(t).
$$
Up to replacing $c_k,d_k$ by $c_ke^{-ik\beta},d_ke^{-ik\beta}$ we may assume that
$\beta=0$ and that $\alpha$ satisfies $e^{ij\alpha}\notin\R$ for $j=1,\ldots,N$.
Then $|v|^2=|u|^2$ and $|\ff_\alpha[v]|^2=|\ff_\alpha[u]|^2$ is equivalent to
\begin{equation}
\label{eq:herm}
\left\{\begin{matrix}
\dst\sum_{j,k=1}^Md_j\overline{d_k}H_j(t)H_k(t)&=&\dst\sum_{j,k=1}^Nc_j\overline{c_k}H_j(t)H_k(t)\\
\dst\sum_{j,k=1}^Md_j\overline{d_k}e^{i(k-j)\alpha}H_j(t)H_k(t)
&=&\dst\sum_{j,k=1}^Nc_j\overline{c_k}e^{i(k-j)\alpha}H_j(t)H_k(t)\\
\end{matrix}\right.
\end{equation}
Looking at the highest order term in \eqref{eq:herm} we obtain $|d_M|^2H_M^2(t)$
on the left hand side and $|c_N|^2H_N(t)^2$ so that $M=N$ and
$|d_M|=|c_N|$.
Up to replacing $v$ by $\frac{c_N}{d_N}v$ we may thus assume that
$d_N=c_N$.

Let us now look at the term of degree $2N-1$ in \eqref{eq:herm}. They appear only in
$H_{N-1}H_N$ thus
$$
\mbox{Re}(d_N\overline{d_{N-1}})=\mbox{Re}(c_N\overline{c_{N-1}})
\quad\mbox{and}\quad
\mbox{Re}(e^{i\alpha}d_N\overline{d_{N-1}})=\mbox{Re}(e^{i\alpha}c_N\overline{c_{N-1}}).
$$
By assumption, $e^{i\alpha}\notin\R$ so that $d_N\overline{d_{N-1}}=c_N\overline{c_{N-1}}$
thus $d_{N-1}=c_{N-1}$.

Let us now assume that $d_{N-j}=c_{N-j}$ for $j=0,\ldots, k-1$ and let us determine
$d_{N-k}$. For this, note that the highest order term
in which $c_{N-k},d_{N-k}$ appear in \eqref{eq:herm} is $H_{N-k}H_N$ which is of order
$2N-k$. As all terms of higher order have already been identified,
the comparison of terms of order $2N-k$ in \eqref{eq:herm} reduces to
$$
\mbox{Re}(d_N\overline{d_{N-k}})=\mbox{Re}(c_N\overline{c_{N-k}})
\quad\mbox{and}\quad
\mbox{Re}(e^{ik\alpha}d_N\overline{d_{N-k}})=\mbox{Re}(e^{ik\alpha}c_N\overline{c_{N-k}}).
$$
By assumption, $e^{ik\alpha}\notin\R$ so that $d_N\overline{d_{N-k}}=c_N\overline{c_{N-k}}$
thus $d_{N-k}=c_{N-k}$.
\end{proof}

\begin{remark}\label{rem:hardy}
One may extend this theorem so that there is no assumption on $v$ by adding a third measure.
More precisely, let $u$ be a Hermite function and assume that
$|\ff_\alpha v|=|\ff_\alpha u|$,
$|\ff_\beta v|=|\ff_\beta u|$ and
$|\ff_\eta v|=|\ff_\eta u|$ for three real numbers $\alpha<\beta<\eta\in[0,\pi]$. Assume that
that $0<\beta-\alpha\leq\pi/2$, $\pi/2\leq\eta-\alpha<\pi$
while $0<\eta-\beta\leq\pi/2$, then $v$ is also a Hermite
function. 

This is a generalised version of Hardy's Uncertainty Principle
which follows immediately from B. Demange's proof of it \cite{De}.
The key point is that the angles $\alpha,\beta,\eta$
define three lines and that none of the angular sectors that they
delimit has an opening of more than $\pi/2$.
Similar results where considered in \cite{HL}.
\end{remark}

\begin{remark}
The hypothesis on $\gamma$ is necessary. For instance, one easily checks that,
if $e^{ij\beta}\in\R$, then $u=h_N+c_{N-j}h_{N-j}$ and $v=h_N+\overline{c_{N-j}}h_{N-j}$
satisfy $|u|=|v|$ and $|\ff_\beta u|=|\ff_\beta v|$.

The proof of Theorem \ref{th:hermite} may easily be adapted to prove that, under the same hypothesis,
if $|\ff_\alpha v|=|\ff_\alpha u|$ and $|\ff_\beta v-\ff_\alpha v|=|\ff_\beta u-\ff_\alpha u|$,
then there exists $c\in\C$ with $|c|=1$ such that $v=cu$.

The second condition may be rewritten $\abs{\frac{\ff_\beta v-\ff_\alpha v}{\beta-\alpha}}
=\abs{\frac{\ff_\beta u-\ff_\alpha u}{\beta-\alpha}}$. It would be tempting to conjecture that the result stays true if we let $\beta\to\alpha$ in the second hypothesis, that is, if we assume that
$|\ff_\alpha v|=|\ff_\alpha u|$ and $\Big|\frac{\partial\ff_\theta v}{\partial\theta}\big|_{\theta=\alpha}\Big|=\Big|\left.\frac{\partial\ff_\theta u}{\partial\theta}\right|_{\theta=\alpha}\Big|$.

This is however not the case since, for $\ff_\alpha u=\sum_{j=0}^Nc_je^{-ij\alpha}h_j$
and
$$
\left.\frac{\partial\ff_\theta u}{\partial\theta}\right|_{\theta=\alpha}=-i\sum_{j=0}^Nc_jje^{-ij\alpha}h_j.
$$
It is then easy to check that $u=h_N+c_{N-1}h_{N-1}$ and $v=h_N+\overline{c_{N-1}}h_{N-1}$
satisfy the hypothesis, but $v$ is not a multiple of $u$ up to a constant phase factor.

We refer to \cite{McDo} for more results in this direction.
\end{remark}

Note that this result may also fall in the scope of Proposition \ref{prop:window} provided one extends it
to distributions since a Hermite function is the convolution of a linear combination of derivatives of Dirac masses at $0$ with a Gaussian.

\begin{remark}
In \cite{BGJ}, we defined a \emph{trivial solution} of the Phase Retrieval Problem as being a linear or antilinear operator that sends $u$ into a solution of the problem. In the problem considered here, this would be a linear or antilinear continuous operator $T$ on $L^2(\R^d)$ such that, for every $u\in L^2(\R^d)$ and every $\alpha\in\tau$, $|\ff_\alpha[Tu]|=|\ff_\alpha u|$. Using the density of Hermite functions, it is not hard to adapt the proof
of \cite[Proposition 3.1]{BGJ} to show that necessarily
$Tu=cu$ for some $c\in\C$, $|c|=1$ as soon as $\tau$ contains at least three numbers $\alpha,\beta,\eta$
satisfying the conditions of Remark \ref{rem:hardy}
\end{remark}

\subsection{Reconstruction of pulse train signals}\ \\
ILet us now consider so called \emph{pulse train} signals which are commonly used in signal processing:

\begin{definition}
A signal $u\in L^2(\R)$ is a (rectangular) \emph{pulse-train signal} (for $H$) if there exists a finite sequence $(a_k)_{k\in\Z}$
such that
$$
u(t)=\sum_{k\in\Z}a_k\chi_{[0,b)}(t-ak)=\sum_{k\in\Z}a_k\chi_{[ak,ak+b)}.
$$
for some $b<a/2$.
\end{definition}

We can now prove the following:

\begin{theorem}\label{th:pulse}
Let $a>0$ and $b<a/2$. Let $(a_k),(b_k)$ be two finite sequences and let
$$
u(t)=\sum_{k\in\Z}a_k\chi_{[ak,ak+b)}\quad\mbox{and}\quad
v(t)=\sum_{k\in\Z}b_k\chi_{[ak,ak+b)}.
$$
Let $\alpha\in\R\setminus\frac{\pi}{2}\Z$.
If $|\ff_\alpha v|=|\ff_\alpha u|$, then there exists a constant $c\in\C$ with $|c|=1$ such that $v=cu$.
\end{theorem}

\begin{remark}
This result is in strong contrast with what happens when $\alpha\in\frac{\pi}{2}\Z$.
In this case, for every $k$, there even exists $u$ such that $|v|=|u|$ and $|\ff v|=|\ff u|$
has at least $k$ solutions, none a constant multiple of the other. We refer to
\cite{Is,Ja} for the construction.
\end{remark}

\begin{remark}
The result should also not be misinterpreted. Theorem \ref{th:pulse} states that a pulse-type signal can be reconstructed
from the modulus of one fractional Fourier transform \emph{among} pulse type signals. In general,
the phase retrieval problem $|\ff_\alpha v|=|\ff_\alpha u|$, $u$ a pulse type signal, will have many other solutions,
see Section \ref{sec:phafr} for a description of these solutions.
\end{remark}

\begin{proof}
There is no loss of generality in assuming that $a=1$ so that $b<1/2$.
An easy computation ({\it see e.g.} \cite{Ja})
then shows that, for $\dst x\in\ent{j-\frac{1}{2},j+\frac{1}{2}}$,
$$
A(u)(x,y)=e^{i\pi j y}\left(\sum_{k\in\Z}a_k\overline{a_{k-j}}e^{2i\pi ky}\right)A(\chi_{[0,b)})(x-j,y)
$$
and a similar expression holds for $A(v)$.
According to \eqref{eq:key}, the hypothesis of the theorem translates into
\begin{eqnarray}
\left(\sum_{k\in\Z}a_k\overline{a_{k-j}}e^{2i\pi ky\cos\alpha}\right)A(\chi_{[0,b)})(-y\sin\alpha-j,y\cos\alpha)&&\notag\\
&&\hspace{-7cm}=\left(\sum_{k\in\Z}b_k\overline{b_{k-j}}e^{2i\pi ky\cos\alpha}\right)A(\chi_{[0,b)})(-y\sin\alpha-j,y\cos\alpha)
\label{eq:pulse}
\end{eqnarray}
for all $j\in\Z$ and for all $y$ such that $-y\sin\alpha\in \ent{j-\frac{1}{2},j+\frac{1}{2}}$
{\it i.e.} for all $y\in I_{\alpha,j}:=\dst\ent{\frac{-2j+1}{2\sin\alpha},\frac{-2j-1}{2\sin\alpha}}$.
Recall that we assumed that $\alpha$ is not a multiple of $\pi$ so that this is perfectly defined.

But, for $|x|\leq b$,
$$
A(\chi_{[0,b)})(x,y)=\dst\frac{e^{i\pi by}}{\pi y}\sin\pi(b-|x|),
$$
thus $A(\chi_{[0,b)})(-y\sin\alpha-j,y\cos\alpha)$ does not vanish on a set $I_j$ of positive measure.
On $I_j$, \eqref{eq:pulse} than reduces to
\begin{equation}
\label{eq:trigpol}
\sum_{k\in\Z}b_k\overline{b_{k-j}}e^{2i\pi ky\cos\alpha}=\sum_{k\in\Z}a_k\overline{a_{k-j}}e^{2i\pi ky\cos\alpha}.
\end{equation}
But, as $(a_k),(b_k)$ are of finite support and as $\cos\alpha\not=0$, this is an equality between
two \emph{trigonometric polynomials}. As this identity is valid on a set $I_j$ of positive measure,
it is valid everywhere:
$$
\sum_{k\in\Z}b_k\overline{b_{k-j}}e^{2i\pi kt}=\sum_{k\in\Z}a_k\overline{a_{k-j}}e^{2i\pi kt}
\mbox{ for every $t\in\R$ and every $j\in\Z$.}
$$
One then easily checks that there is $c\in\C$
with $|c|=1$ such that $b_k=ca_k$ for every $k$, thus $v=cu$.
\end{proof}

\begin{remark}
Rectangular pulse trains may be replaced by more general pulse trains
$u(t)=\sum_{k\in\Z}a_kH(t-ak)$ where $H$ is supported in $[0,b)$ $b<a/2$, provided that,
for every $j\in\Z$, $A(H)(-y\sin\alpha-j,y\cos\alpha)$ does not vanish on a set $I_j$ of positive measure.
This is the case if $H(t)=e^{i(\alpha t+\beta t^2)}\chi_{[0,b)}(t)$, $\alpha,\beta\in\R$ since then
$A(H)(x,y)=e^{i\alpha x+i\beta x^2/2}A(\chi_{[0,b)})(x,y-\beta x/\pi)$.

Further, as $A(H)$ is continuous and $A(H)(0,y)=\ff[|H|^2](y)$, it is enough to have $$
\ff[|H|^2](j\tan\alpha)\not=0\mbox{ for }j\in\Z.
$$
This is easily seen to be the case when $H=(1-t/b)_+$.
\end{remark}

\subsection{Reconstruction of linear combinations of translates of Gaussians}
In this section, we consider signals of the form:

\begin{definition}
A function $u\in L^2$ is a \emph{combination of shifted Gaussians}
if there exist an integer $N$, complex numbers $c_1,\ldots,c_N$ and real numbers
$t_1,\ldots,t_N$ such that
$$
u(t)=\sum_{j=1}^Nc_j\gamma(t-t_j)
$$
where $\gamma(t)=e^{-\pi t^2}$.
\end{definition}

Note that, if $\delta_a$ is the dirac mass at $a$, then $u=\dst\left(\sum_{j=1}^Nc_j\delta_{t_j}\right)*\gamma$. Thus, the next theorem has to be compared to Proposition \ref{prop:window} (extended to distributions).

\begin{theorem}\label{th:510}
Let $\alpha\in\R\setminus\dst\frac{\pi}{2}\Z$ and let $u,v$ be two combination of shifted Gaussians.
If $|\ff_\alpha v|=|\ff_\alpha u|$, there exists $c\in\C$ with $|c|=1$ such that $v=cu$.
\end{theorem}

\begin{proof}
Let us write
$$
u(t)=\sum_{j=1}^Nc_j\gamma(t-t_j)\quad\mbox{and}\quad
v(t)=\sum_{j=1}^M\kappa_j\gamma(t-\tau_j)
$$
with $c_j,\kappa_j\in\C$ and $t_j,\tau_j\in\R$. Moreover the $t_j$'s (resp. the $\tau_j$'s)
are all distinct.

An easy computation shows that $A(\gamma)(x,y)=2^{-1/2}e^{-\pi(x^2+y^2)/2}$, so that
with Property \eqref{prop:amb3b} of the ambiguity function and its bilinearity, we get
$$
A(u)(x,y)=\frac{1}{\sqrt{2}}\sum_{j,k=1}^Nc_j\overline{c_k}e^{-\pi(t_j+t_k)^2/2}e^{-\pi\bigl((x+t_k-t_j)^2+(y-it_k-it_j)^2\bigr)/2}.
$$
In particular,
\begin{eqnarray}
A(u)(-t\sin\alpha,t\cos\alpha)
&=&\frac{1}{\sqrt{2}}\sum_{j,k=1}^Nc_j\overline{c_k}e^{-\pi(t_j+t_k)^2/2}e^{-\pi\bigl((-t\sin\alpha+t_k-t_j)^2+(t\cos\alpha+it_k+it_j)^2\bigr)/2}\nonumber\\
&=&\frac{e^{-\pi t^2/2}}{\sqrt{2}}\sum_{j,k=1}^Nc_j\overline{c_k}e^{-\pi(t_j-t_k)^2/2}e^{-i\pi t\bigl(t_je^{i\alpha}+t_ke^{-i\alpha}\bigr)}.
\label{eq:sgausu}
\end{eqnarray}
A similar expression holds for $v$:
\begin{equation}
\label{eq:sgausv}
A(v)(-t\sin\alpha,t\cos\alpha)=2^{-1/2}
e^{-\pi t^2/2}\sum_{j,k=1}^M\kappa_j\overline{\kappa_k}e^{-\pi(\tau_j-\tau_k)^2/2}
e^{-i\pi t\bigl(\tau_je^{i\alpha}+\tau_ke^{-i\alpha}\bigr)}.
\end{equation}

We will now need the two following facts:

\medskip

\noindent{\bf Fact 1.} {\sl Let $\alpha\in\R\setminus\dst\frac{\pi}{2}\Z$ and let
$\{t_j\}_{j\in\Z}$ and $\{\tau_j\}_{j\in\Z}$ be two finite sequences of real numbers. If there exists
$j,k,j',k'\in\Z$ such that 
$t_je^{i\alpha}+t_ke^{-i\alpha}=\tau_{j'}e^{i\alpha}+\tau_{k'}e^{-i\alpha}$
then $t_j=\tau_{j'}$ and $t_k=\tau_{k'}$.} 

\medskip

\begin{proof}[Proof of Fact 1]
The condition is equivalent to
\begin{equation}
\label{eq:simple}
(t_j-\tau_{j'})e^{2i\alpha}=\tau_{k'}-t_k.
\end{equation}
As $\tau_{k'}-t_k\in\R$, \eqref{eq:simple} implies that $(t_j-\tau_{j'})e^{2i\alpha}\in\R$.
But, we assumed that $\alpha\in\R\setminus\dst\frac{\pi}{2}\Z$, so that $e^{2i\alpha}\notin\R$,
thus $t_j-\tau_{j'}=\tau_{k'}-t_k=0$.
\end{proof}

\medskip

\noindent{\bf Fact 2.} {\sl The set $\{e^{zt}\}_{z\in\C}$ is 
linearly independent set of functions on $\R$.} 

\medskip

\begin{proof}[Proof of Fact 2]
Let us consider a finite linear combination of $e^{zt}$ that vanishes:
$$
G(t):=\sum_{j=1}^N\lambda_je^{z_j t}=0\quad\lambda_j,z_j\in\C
$$
and the $z_j$'s are all distinct. 
Then evaluating the derivatives $G,G',\ldots,G^{(N-1)}$ at $0$, we obtain the Vandermonde system
$$
\left\{\begin{matrix}\lambda_1&+&\lambda_2&+&\cdots&+&\lambda_N&=&0\\
\lambda_1z_1&+&\lambda_2z_2&+&\cdots&+&\lambda_Nz_N&=&0\\
\vdots&&\vdots&&\ddots&&\vdots&&\\
\lambda_1z_1^{N-1}&+&\lambda_2z_2^{N-1}&+&\cdots&+&\lambda_Nz_N^{N-1}&=&0\\
\end{matrix}\right..
$$
As the $z_j$'s are all distinct, this system has non zero determinant, thus $\lambda_1=\cdots=\lambda_N=0$. 
\end{proof}

We can now complete the proof of the theorem.
From Fact 1, each term of the form $e^{zt}$
appearing in the sum \eqref{eq:sgausu} --resp. \eqref{eq:sgausv}-- appears exactly ones.
Moreover, as $|\ff_\alpha v|=|\ff_\alpha u|$
implies $A(v)(-t\sin\alpha,t\cos\alpha)=A(u)(-t\sin\alpha,t\cos\alpha)$,
Fact 2 implies that the two sums are equal term by term:
$M=N$, $\{\tau_j\}=\{t_j\}$ thus (up to renumbering) $\tau_j=t_j$
and then, for every $j,k=1,\ldots,N$, $\kappa_j\overline{\kappa_k}=c_j\overline{c_k}$.
This last identity implies that there exists $c\in\C$ with $|c|=1$ such that
$\kappa_j=cc_j$ for all $j$ and finally that $v=cu$.
\end{proof}

\begin{remark}
Again the result should not be over-interpreted and it only says that uniqueness is achieved
from the measure of $|\ff_\alpha u|$ \emph{among} combinations of shifted Gaussians,
provided $\alpha\notin\dst\frac{\pi}{2}\Z$.
Again this is in strong contrast with the Pauli problem for which this result is false.
\end{remark}

We have only used the fact that the shifts $t_j$'s and $\tau_j$'s are real in a mild fashion.
Actually the same proof works if the $t_j$'s and $\tau_j$'s are all purely imaginary or even 
if we restrict them to be in a set $E\in\C$ such that its sum-set
$E+E=\{e+e'\,:\ e,e'\in E\}$ intersects its $2\alpha$-rotate only at $0$:
$\bigl(e^{2i\alpha}(E+E)\bigr)\cap(E+E)=\{0\}$. Of course, taking $t_j$ or $\tau_j$ to be complex amounts
to taking time-frequency shifts of Gaussians. Let us now prove a more general theorem in which we further 
also allow for dilates:

\begin{theorem}
Let $\alpha\in]0,\frac{\pi}{2}[$ and let 
$u,v$ be time-frequency translates of Gaussians, that is, $$
u(t)=\sum_{j=1}^N \sum_{k=1}^{N_j}c_{j,k}e^{2i\pi\omega_ {j,k} t}\gamma\left(\frac{t-s_{j,k}}{\sigma_j^{1/2}}\right)
\quad\mbox{and}\quad
v(t)=\sum_{j=1}^M  \sum_{k=1}^{M_j}\kappa_{j,k}e^{2i\pi\eta_{j,k} t}\gamma\left(\frac{t-t_{j,k}}{\tau_j^{1/2}}\right)
$$
with $c_{j,k},\kappa_{j,k}\in\C$, $s_{j,k},t_{j,k},\omega_{j,k},\eta_{j,k}\in\R$,
and $\sigma_1>\cdots>\sigma_N>0$, $\tau_1>\cdots>\tau_M>0$.

Let us define the following set of angles: $\dst\alpha_0^\pm=\pm\left(\frac{\pi}{2}-\alpha\right)$;
and for each $j$ for which there exists $k$ such that $\omega_{j,k}\not=0$ (resp. $\eta_{j,k}\not=0$),
let $\alpha_j$ (resp. $\gamma_j$) be defined by
$$
e^{i\alpha_j}=\frac{\sigma_j^2\cos\alpha+i\sin\alpha}{\sqrt{\alpha_j^4+1}}\quad\mbox{resp.}\quad
e^{i\gamma_j}=\frac{\tau_j^2\cos\alpha+i\sin\alpha}{\sqrt{\tau_j^4+1}}.
$$
Assume that $\alpha_0^\pm,\alpha_1\ldots,\alpha_N$ are all distinct and that
$\alpha_0^\pm,\gamma_1\ldots,\gamma_N$ are all distinct.

Then, $|\ff_\alpha v|=|\ff_\alpha u|$
implies that there exists $c\in\C$ with $|c|=1$ such that $v=cu$.
\end{theorem}

\begin{proof} The proof follows in part the lines of the previous theorem. 
For simplicity of notation, let us write
$$
u(t)=\sum_{j=1}^N c_je^{2i\pi\omega_j t}\gamma\left(\frac{t-s_j}{\sigma_j^{1/2}}\right)
\quad\mbox{and}\quad
v(t)=\sum_{j=1}^M \kappa_je^{2i\pi\eta_j t}\gamma\left(\frac{t-t_j}{\tau_j^{1/2}}\right)
$$
where $c_j,\kappa_j\in\C$, $\sigma_j,\tau_j>0$ and $s_j,t_j,\omega_j,\eta_j\in\R$. 
As previously, $|\ff_\alpha v|=|\ff_\alpha u|$ implies that $A(v)(-t\sin\alpha,t\cos\alpha)=A(u)(-t\sin\alpha,t\cos\alpha)$. We will thus compute
$A(u)$ and $A(v)$ and compare the expressions obtained.

\medskip

Let us first note that, writing $\gamma_{\sigma,a,\omega}(t)=e^{2i\pi\omega t}\gamma\left(\frac{t-a}{\sigma^{1/2}}\right)$, we obtain
\begin{eqnarray*}
&&
A(\gamma_{\sigma,a,\omega},\gamma_{\tau,b,\eta})(x,y)=\left(\frac{\sigma\tau}{\sigma+\tau}\right)^{1/2}
e^{i\pi\bigl((\omega+\eta)x+(a+b)(y-\omega+\eta)\bigr)}\\
&&\times\exp-\frac{\pi}{\sigma+\tau}\left((x-a+b)^2+\sigma\tau (y-\omega+\eta)^2-i(\tau-\sigma)(x-a+b)(y-\omega+\eta)\right).
\end{eqnarray*}
In particular, $A(\gamma_{\sigma,a,\omega},\gamma_{\tau,b,\eta})(-t\sin\alpha,t\cos\alpha)$ is expressed in the form $c\exp -P_{\sigma,a,\omega,\tau,b,\eta}(t)$ where
$P_{\sigma,a,\omega,\tau,b,\eta}(t)$ is a polynomial of degree $2$ with highest order term
$$
\pi \frac{\sin^2\alpha+\sigma\tau\cos^2\alpha-i(\sigma-\tau)\sin\alpha\cos\alpha}{\sigma+\tau}t^2.
$$

Let us denote by $\mathbb{C}_2[X]$, the set of complex polynomials of degree $2$. For $P\in\mathbb{C}_2[X]$, let us write
$J_P(u)=\{(j,k)\,:P_{\sigma_j,s_j,\omega_j,\sigma_k,s_k,\omega_k}=P\}$ and
$J_P(v)=\{(j,k)\,:P_{\tau_j,t_j,\eta_j,\tau_k,t_k,\eta_k}=P\}$. We may then write
$$
A(u)(-t\sin\alpha,t\cos\alpha)=\sum_{P\in\mathbb{P}_2[\C]}
\sum_{(j,k)\in J_P(u)}\left(\frac{\sigma_j\sigma_k}{\sigma_j+\sigma_k}\right)^{1/2}c_j\overline{c_k}e^{-P(t)}.
$$
A similar expression holds for $A(v)$ with $J_P(v)$ instead of $J_P(u)$.
We will now need an elaboration on Fact 2:

\medskip

\noindent{\bf Fact 3.} {\sl The set of functions $\{e^{zt+\zeta t^2}\}_{z,\zeta\in\C}$
is a set of linearly independent functions over $\R$.}

\medskip

We postpone the proof of this fact to the end of the proof of the theorem.
As a consequence of this fact, we obtain that, for each $P\in\mathbb{C}_2[X]$,
$$
\sum_{j,k\in J_P(u)}\left(\frac{\sigma_j\sigma_k}{\sigma_j+\sigma_k}\right)^{1/2}c_j\overline{c_k}e^{-P(t)}
=\sum_{j,k\in J_P(v)}\left(\frac{\tau_j\tau_k}{\tau_j+\tau_k}\right)^{1/2}\kappa_j\overline{\kappa_k}e^{-P(t)}.
$$

The following fact will allow us to get some information on the sets $J_P(u)$ and $J_P(v)$\,:

\medskip

\noindent{\bf Fact 4.} {\sl Let $\alpha\in]0,\pi/2[$ and assume that
$u,v,u',v'>0$ are such that 
\begin{eqnarray}
&&\frac{\sin^2\alpha+uv\cos^2\alpha-i(u-v)\sin\alpha\cos\alpha}{u+v}\nonumber\\
&&\qquad\qquad\qquad=\frac{\sin^2\alpha+u'v'\cos^2\alpha-i(u'-v')\sin\alpha\cos\alpha}{u'+v'}.
\label{eq:Fact4}
\end{eqnarray}
Then $u'=u$ and $v'=v$.}

\medskip

\begin{proof}[Proof of Fact 4] Note that the condition on $\alpha$ is simply $\cos\alpha\not=0$ and $\sin\alpha\not=0$.

Looking at the imaginary part in \eqref{eq:Fact4} gives
$\dst\frac{u-v}{u+v}=\frac{u'-v'}{u'+v'}$, that is $\frac{1-v/u}{1+v/u}=\frac{1-v'/u'}{1+v'/u'}$.
But, as one easily checks, $\dst t\to\frac{1-t}{1+t}$ is (strictly) decreasing on $[0,+\infty)$, thus one-to-one, thus $v'/u'=v/u$.

Now, looking at the real part in \eqref{eq:Fact4} gives
$$
\frac{\sin^2\alpha+uv\cos^2\alpha}{u+v}
=\frac{\sin^2\alpha+u'v'\cos^2\alpha}{u'+v'}.
$$
Factoring out $u$ and $u'$, we obtain
$$
\frac{\frac{\sin^2\alpha}{u}+u/v\cos^2\alpha}{1+v/u}
=\frac{\frac{\sin^2\alpha}{u'}+u'/v'\cos^2\alpha}{1+v'/u'}.
$$
As $v'/u'=v/u$ and $\sin\alpha\not=0$, we get $1/u'=1/u$ and then $v'=v$.
\end{proof}

\medskip

As a consequence, if $(j,k),(j',k')\in J_P(u)$ and $(j'',k'')\in J_P(v)$, then 
$\sigma_j=\sigma_{j'}=\tau_{j''}$ and $\sigma_k=\sigma_{k'}=\tau_{k''}$.
We may therefore group all terms stemming from a given $\sigma$ in the expansions of
$A(u)$ and $A(v)$ and identify those terms.

In other words, fix $\sigma>0$, and define
$$
u_\sigma(t)=\sum_{\begin{matrix}\scriptstyle j=1,\ldots,N\\ \scriptstyle \sigma_j=\sigma\end{matrix}} c_je^{2i\pi\omega_j t}\gamma\left(\frac{t-s_j}{\sigma^{1/2}}\right)
$$
and
$$
v_\sigma(t)=\sum_{\begin{matrix} \scriptstyle j=1,\ldots,M\\ \scriptstyle \tau_j=\sigma\end{matrix}} \kappa_je^{2i\pi\eta_j t}\gamma\left(\frac{t-t_j}{\sigma^{1/2}}\right).
$$
Then
$$
A(u)(-t\sin\alpha,t\cos\alpha)=A(u_\sigma)(-t\sin\alpha,t\cos\alpha)+R_u(t)
$$
and
$$
A(v)(-t\sin\alpha,t\cos\alpha)=A(v_\sigma)(-t\sin\alpha,t\cos\alpha)+R_v(t)
$$
where $R_u$ and $R_v$ are each linearly independent both from $A(u_\sigma)(-t\sin\alpha,t\cos\alpha)$
and from $A(v_\sigma)(-t\sin\alpha,t\cos\alpha)$. As
$A(u)(-t\sin\alpha,t\cos\alpha)=A(v)(-t\sin\alpha,t\cos\alpha)$, this implies that
$A(u_\sigma)(-t\sin\alpha,t\cos\alpha)=A(v_\sigma)(-t\sin\alpha,t\cos\alpha)$.
In other words, we may now assume that $u=u_\sigma$ and $v=v_\sigma$ or, equivalently, that
$\sigma_j=\tau_{j'}:=\sigma$ for every $j,j'$. But then
\begin{equation}
\label{eq:au2}
A(u)(-t\sin\alpha,t\cos\alpha)
=\sqrt{\frac{\sigma}{2}}\sum_{j,k=1}^Nc_j\overline{c_k}e^{-P_{\sigma,s_j,\omega_j,\sigma,s_k,\omega_k}(t)}
\end{equation}
and
\begin{equation}
\label{eq:av2}
A(v)(-t\sin\alpha,t\cos\alpha)
=\sqrt{\frac{\sigma}{2}}\sum_{j,k=1}^M\kappa_j\overline{\kappa_k}e^{-P_{\sigma,t_j,\eta_j,\sigma,t_k,\eta_k}(t)}.
\end{equation}

We have already exploited the fact that the $P_{\sigma,s_j,\omega_j,\sigma,s_k,\omega_k}$'s
and $P_{\sigma,t_j,\eta_j,\sigma,t_k,\eta_k}$'s all have same second order term. Let us now use again Fact 3 (which in this case simplifies to Fact 2) and compare those terms for which those polynomials have same first order term.

Note that the first order term $P_{\sigma,a,\omega,\sigma,b,\eta}(t)$ is
$$
\pi\left(ae^{i(\pi/2+\alpha)}+be^{i(\pi/2-\alpha)}
-\eta\sqrt{\sigma^4+1}e^{i\alpha_\sigma}+\omega\sqrt{\sigma^4+1}e^{-i\alpha_\sigma}\right)t
$$
where $e^{i\alpha_\sigma}=\dst\frac{\sigma^2\cos\alpha+i\sin\alpha}{\sqrt{\sigma^4+1}}$.
But, as in the previous proof, our assumption on the angles implies that the two sums \eqref{eq:au2}
and \eqref{eq:av2} are equal term by term. In other words $M=N$,
$(s_j,\omega_j)=(t_j,\eta_j)$ for all $j$'s (up to reordering) and $c_j\overline{c_k}=\kappa_j\overline{\kappa_k}$ for every $j,k$,
thus $c_j=c\kappa_j$ with $|c|=1$.
\end{proof}

Before proving Fact 3, we need an intermediate elaboration on Fact 2:

\medskip

\begin{proof}[Proof of Fact 3] 
Let us assume that there is a vanishing linear combination with non zero coefficients $\mu_k$:
\begin{equation}
\label{eq:linind1}
\sum_{k=1}^N\mu_ke^{-\zeta_k t^2+z_kt}=0\quad\mbox{for }t\in\R.
\end{equation}
In the previous notation, we of course assume that no term is repeated $(\zeta_k,z_k)\not=(\zeta_l,z_l)$.
As this is an equality of entire functions, it is true for all $t\in\C$.

Without loss of generality, we assume that $|\zeta_1|=\max |\zeta_k|$ and that
$\zeta_1$ is real non negative --- by changing $t$ to $te^{-i(\arg\zeta_1)/2}$.
We also reorder the remaining terms to have $\zeta_1=\cdots=\zeta_M$ while all other $\zeta_k$'s are
$\not=\zeta_1$. We then re-order the $M$ first terms to have $\Re(z_1)=\cdots=\Re(z_L)=\max\Re(z_k)$
and $\Re(z_k)<\Re(z_1)$ for $k>L$.
Further re-order the $L$ first terms to have $\Im(z_1)<\Im(z_2),\ldots,\Im(z_L)$.
Note that this implies that if $\alpha>0$ is small enough, then $\Re(z_1e^{i\alpha/2})>\Re(z_ke^{i\alpha/2})$
for $k=1,\ldots,M$.\footnote{This means that we choose those $k$'s for which $\zeta_k$ has largest modulus,
rotate so that $\zeta_k$ has also largest real part.
Then we chose among them, those $k$'s for which $\Re(z_k)$ is the largest and then, after a small rotation,
such that $z_1$ is the unique point for which $\Re(z_k)$ is the largest.}

Notice that $\zeta_1=\max |\zeta_k|$ implies that, for $k\geq M+1$, $\Re(\zeta_1-\zeta_k)>0$.
Therefore, there is an $\alpha_0$ such that, for $|\alpha|<\alpha_0$ and $k\geq M+1$, $\Re\bigl((\zeta_1-\zeta_k)e^{i\alpha}\bigr)>0$. It follows that
$e^{-(\zeta_k-\zeta_1) (e^{i\alpha/2}t)^2+(z_k-z_1)e^{i\alpha/2}t} \to 0$ as $t\to\pm\infty$.

Factoring out $e^{-\zeta_1 t^2+z_1t}$ in \eqref{eq:linind1} and setting $t=e^{i\alpha/2}\tau$, we obtain
$$
e^{-\zeta_1 e^{i\alpha} \tau^2+z_1e^{i\alpha/2}\tau}\sum_{k=1}^N\mu_ke^{-(\zeta_k-\zeta_1)e^{i\alpha} \tau^2+(z_k-z_1)e^{i\alpha/2}\tau}=0.
$$
Dividing by $e^{-\zeta_1 e^{i\alpha} \tau^2+z_1e^{i\alpha/2}\tau}$ and letting $\tau\to\pm\infty$, we obtain
\begin{equation}
\label{eq:linind2}
\sum_{k=1}^M\mu_ke^{(z_k-z_1)e^{i\alpha/2}\tau}
=\mu_1+\sum_{k=2}^M\mu_ke^{(z_k-z_1)e^{i\alpha/2}\tau}\to 0\mbox{ when }\tau\to\pm\infty\quad\mbox{for }|\alpha|\leq\alpha_0.
\end{equation}
As explained above, by taking $\alpha>0$ small enough, $\Re\big((z_k-z_1)e^{i\alpha/2}\big)<0$ if $k=2,\ldots,M$,
thus \eqref{eq:linind2} takes the form $\mu_1+\ffi(\tau)\to0$ where $\ffi(\tau)\to0$ as $\tau\to+\infty$.
Therefore $\mu_1=0$, a contradiction.
\end{proof}

\section{Conclusion}

In this paper, we have studied the phase retrieval problem for single or multiple measurements of the
fractional Fourier transform. This problem occurs naturally in quantum mechanics and in optics.

For the single measurement case and compactly supported functions, we have shown that the zero-flipping phenomena occurs as for the usual Fourier transform. However, if one is interested in the problem
for more structured signals like translates of Gaussians, then uniqueness is achieved.

For multiple measurements, we have seen that uniqueness is achieved in many cases:
for Hermite signals and pulse train signals, two measurements suffice, provided the parameters of the FrFT are chosen properly. For compactly supported functions, a countable set of measurements guaranties uniqueness and a reconstruction Formula is provided.

In a forthcoming paper, we will explore the practical validity of this reconstruction formula.
We will also propose a modification of the Grechter-Saxton iterative algorithm
({\it see e.g.} \cite{Fi}) for multiple phase-less FrFT measurements and study its validity.

%

\end{document}